\pdfoutput=1
\documentclass[a4paper,11pt]{amsart}
\usepackage[english]{babel}
\usepackage{amsfonts}
\usepackage{amsthm}
\usepackage{amsmath}
\usepackage{amssymb}
\usepackage{enumerate}
\usepackage{enumitem}
\usepackage{mathtools}
\usepackage{tikz-cd}
\usepackage[all]{xy}
\usepackage{graphicx}
\usepackage{amsmath}
\usepackage{extpfeil}
\usepackage{comment}
\usepackage{float}
\usepackage[labelformat=empty]{caption}
\usepackage[toc]{appendix}
\usepackage[left=3cm,top=2.5cm,right=3cm,bottom=2.5cm]{geometry}
\usepackage{hyperref}
\usepackage{comment}
\usepackage{resizegather}
\usepackage[activate={true, nocompatibility}, final, tracking=true, kerning=true, spacing=true]{microtype}
\usepackage{xcolor}
\usepackage{bm}

\theoremstyle{plain}
\newtheorem{thm}{Theorem}[section]
\newtheorem{coroll}[thm]{Corollary}
\newtheorem{lemma}[thm]{Lemma}

\newtheorem{example}[thm]{Example}
\newtheorem{prop}[thm]{Proposition}
\newtheorem{prin}[thm]{Principle}

\newtheorem{defn}[thm]{Definition}
\newtheorem{remark}[thm]{Remark}

\def\makeCal#1{%
\expandafter\newcommand\csname c#1\endcsname{\mathcal{#1}}}
\def\makeBB#1{%
\expandafter\newcommand\csname b#1\endcsname{\mathbb{#1}}}
\def\makeFrak#1{%
\expandafter\newcommand\csname f#1\endcsname{\mathfrak{#1}}}

\newcounter{int}
\loop
\addtocounter{int}{1}
\edef\y{\Alph{int}}%
\expandafter\makeCal\y
\expandafter\makeBB\y
\expandafter\makeFrak\y
\ifnum\value{int}<26%
\repeat

\makeatletter
\newcommand{\colim@}[2]{%
  \vtop{\m@th\ialign{##\cr
    \hfil$#1\operator@font colim$\hfil\cr
    \noalign{\nointerlineskip\kern1.5\ex@}#2\cr
    \noalign{\nointerlineskip\kern-\ex@}\cr}}%
}
\newcommand{\colim}{%
  \mathop{\mathpalette\colim@{\rightarrowfill@\textstyle}}\nmlimits@
}
\makeatother

\tikzset{
  symbol/.style={
    draw=none,
    every to/.append style={
      edge node={node [sloped, allow upside down, auto=false]{$#1$}}}
  }
}

\newcommand{\dash}{{\operatorname{-}}}

\microtypecontext{spacing=nonfrench}

\usepackage[normalem]{ulem}
\usepackage{pbox}

\DeclareMathOperator{\ad}{ad}
\DeclareMathOperator{\Alg}{Alg}
\DeclareMathOperator{\Aut}{Aut}
\DeclareMathOperator{\Ann}{Ann}

\DeclareMathOperator{\Der}{Der}
\DeclareMathOperator{\diag}{diag}
\newcommand{\dual}{\vee}
\DeclareMathOperator{\End}{End}

\DeclareMathOperator{\GL}{GL}

\DeclareMathOperator{\Hom}{Hom}
\newcommand{\id}{\mathrm{id}}

\DeclareMathOperator{\JRad}{J}

\DeclareMathOperator{\LieRad}{Rad}
\DeclareMathOperator{\LieNilRad}{NilRad}

\DeclareMathOperator{\Mat}{Mat}

\DeclareMathOperator{\semistable}{ss}

\DeclareMathOperator{\Spec}{Spec}
\DeclareMathOperator{\Split}{Split}

\DeclareMathOperator{\Sym}{Sym}

\DeclareMathOperator{\Trace}{Tr}
\DeclareMathOperator{\Transpose}{t}

\DeclareMathOperator{\Gr}{Gr}
\DeclareMathOperator{\GTrace}{GTr}

\DeclareMathOperator{\wt}{wt}

\DeclareMathOperator{\ModAlg}{\Alg(V)}

\DeclareMathOperator{\ModLieAlg}{Lie(V)}
\DeclareMathOperator{\Span}{Span}
\DeclareMathOperator{\Can}{Can}

\DeclareMathOperator{\numinvariant}{\nu}
\DeclareMathOperator{\basefield}{\mathit{k}}
\DeclareMathOperator{\Nik}{Nik}
\DeclareMathOperator{\riesz}{riesz}

\begin{document}
\thispagestyle{empty}

\title{Canonical Filtrations of Finite-Dimensional Algebras}

\begin{abstract}
We study canonical filtrations of finite-dimensional associative algebras and Lie algebras. These filtrations are defined via optimal destabilizing one-parameter subgroups in the sense of geometric invariant theory (GIT), and appear to be a new invariant of finite-dimensional algebras. We establish some fundamental properties of these filtrations, and show that an algebra is semisimple if and only if it is GIT semistable. We give a method to compute canonical filtrations of algebras whose automorphism group or module of derivations is sufficiently rich, and use this to compute these filtrations in examples when the automorphism group contains a sufficiently large torus. We also obtain some results on the structure of the associated graded algebra of the canonical filtration of an associative algebra.
\end{abstract}

\author{Trevor Jones}

\maketitle

\section{Introduction}

One of the goals of moduli theory is to construct spaces which classify all objects of a given type. In order to obtain a space with reasonable properties (e.g. finite-dimensional spaces), typically we must consider a restricted class of semistable objects. To define a notion of semistability, one usually defines a function $\omega_x$ on the set of all filtrations of an object $x$, and then declares $x$ to be semistable if $\omega_x(f) \geq 0$ for every filtration $f$ of $x$. The unstable points are then the objects $x$ for which there is a filtration $f$ such that $\omega_x(f) < 0$. In good circumstances, as is the case when the functions $\omega_x$ come from geometric invariant theory (GIT), every unstable object $x$ has an essentially unique canonical filtration which is ``most destabilizing'' in the sense that it minimizes a normalized version of the function $\omega_x$. It is then natural to ask if there are intrinsic descriptions of the semistable objects and the canonical filtrations of unstable objects. These questions have been the object of much study when the underlying objects form an Abelian category, but much less is known in the non-Abelian setting.

In this work, we study these questions for non-Abelian categories of finite-dimensional algebras and associative algebras over a perfect field $\basefield$. Our goal is to investigate to what extent we can intrinsically characterize (1) the semistable algebras and (2) the canonical filtrations of unstable algebras. A filtration of an algebra $A$ in this setting is the data of a weighted flag of subspaces
\[ A = F_{0} \supsetneq F_{1} \supsetneq \ldots \supsetneq F_{N} \]
where the weights $w_i$ are increasing in $i$ and the filtration is compatible with the multiplication or Lie bracket $\mu_A$ in the sense that $\mu_A(F_i, F_j) \subseteq F_{r}$ implies $w_r \geq w_i + w_j$. Semistability and the canonical filtration of $A$ are characterized via a numerical invariant $\nu(A, \dash)$ (see definition (\ref{defn: numerical invariants}) for details) on the space of filtrations of $A$:
\begin{itemize}
    \item $A$ is semistable if $\nu(A, (F_\bullet, w_\bullet)) \geq 0$ for all non-trivial filtrations $(F_\bullet, w_\bullet)$.
    \item The canonical filtration $F(A)^{\Can}_\bullet$ of $A$ is the trivial filtration if $A$ is semistable and the filtration which minimizes $\nu(A, \dash)$ if $A$ is unstable.
\end{itemize}

There are several difficulties that arise when computing canonical filtrations. First, is not clear how to efficiently search for an optimizer of $\nu$. In Abelian categories there is a standard procedure to iteratively construct canonical filtrations (e.g. Harder-Narasimhan filtrations) from ``optimally destabilizing'' subobjects; see \cite{rudakov1997stability}. This procedure doesn't generalize to the non-Abelian setting since quotients by arbitrary subalgebras need not exist. Second, given an unweighted sequence of subspaces $F_{\bullet}$ of $A$ as above, determining the weights $w_\bullet$ which minimize $\nu(A, \dash)$ on the convex cone determined by $F_\bullet$ requires solving a constrained optimization problem for which there is no closed-form solution. Furthermore, a priori it is not clear how to restrict the set of filtrations one must search over, so computing the canonical filtration requires solving many optimization problems of this form. 

One of our key observations is that the canonical filtration is invariant under automorphisms and derivations. This enables us to give a complete characterization of the semistable algebras, achieving goal (1) and generalizing \cite[Theorem 4.3]{lauret2003moment} and \cite[Theorem 4.6]{zhang2023moment} to perfect fields.:

\begin{thm}[= Theorems \ref{thm: A Lie algebra is semistable iff semisimple}, \ref{thm: an associative algebra is semistable iff semisimple}]
A Lie algebra or associative algebra $A$ is semistable if and only if it is semisimple.
\end{thm}

Intrinsically characterizing the canonical filtration of an unstable algebra appears to be a difficult problem, making goal (2) challenging. In fact, to the best of our knowledge these canonical filtrations do not correspond to any known algebra filtrations (e.g. the filtration given by powers of the Jacobson radical), and thus appear to be new invariants for finite-dimensional algebras. Even for relatively simple algebras whose canonical filtration we can compute, such as algebras of the form $k[x_1, \ldots, x_n] / (x_1, \ldots, x_n)^d$, the filtration can behave in surprising ways. Examples (\ref{example: small monomial algebra whose canonical filtration does not arise from a grading}) and (\ref{example: S_{n,4}}) illustrate some of these phenomena.

Despite the difficulty of computing the canonical filtration in general, when an algebra has a sufficiently rich automorphism group or module of derivations we can place stronger constraints on this filtration. We use this method to great effect for $\bZ^r$-graded algebras and show that the filtered pieces of the canonical filtration are spanned by homogeneous elements. For $\bZ^r$-graded algebras $A$ with the property that each homogeneous component $A_\alpha$ is one dimensional, this implies that we can compute the canonical filtration by solving a single optimization problem. Classes of algebras with such a grading include commutative monomial algebras, algebras of upper triangular matrices, and Borel subalgebras of semisimple Lie algebras over $\bC$. We use this method to compute the canonical filtration in some explicit examples for these types of algebras.

We also obtain results on the structure of graded-semistable associative algebras, partially generalizing \cite[Theorem 4.10]{zhang2023moment} to perfect fields using the invariance of the canonical filtration under automorphism and derivations. We recall that a $\bZ$-graded algebra is \emph{graded-semistable} if its canonical filtration is the split filtration induced by the grading. Such algebras are closely related to canonical filtrations due to the Recognition Theorem (\cite[Theorem 5.4.4]{halpern2014structure}), which in our context states that a filtration $(F_\bullet, w_\bullet)$ is the canonical filtration if and only if the associated graded algebra $\Gr(F_\bullet, w_\bullet)$ is graded-semistable. 

While this paper was in preparation, we became aware of the two papers \cite{lauret2003moment} and \cite{zhang2023moment}, which study Lie algebra and associative algebra structures on a finite-dimensional complex vector space $V$ from the symplectic viewpoint, respectively. Specifically, they study the critical points of the norm-squared of a certain moment map on the projective space $\bP(\Hom(V^{\otimes2}, V))$. While we recover similar results on the structure of graded-semistable associative algebras (which correspond to the critical points mentioned above), the questions we consider are more general than those investigated in these works. We study canonical filtrations for arbitrary algebras, not just those corresponding to critical points, and also work over an arbitrary perfect field instead of the complex numbers. Our methods are more algebraic in nature, and show that the invariance of the canonical filtration under automorphisms and derivations provides a fundamental means to understand these filtrations.

\subsection{Contents}

In \autoref{section: preliminaries}, we set up the GIT problems by defining the parameter schemes $\ModAlg$ and $\ModLieAlg$ of associative and Lie algebra structures on $V$, a $\GL(V)$ action on these schemes such that orbits correspond to isomorphism classes of algebras, the semistability condition (\ref{defn: semistablity for algebras}), the numerical invariant $\numinvariant$ (\ref{defn: numerical invariants}), and the canonical filtration (\ref{defn: the canonical filtration}) of an algebra.
In \autoref{section: properties of the canonical filtration}, we establish some key properties of the canonical filtration, the most important being that the canonical filtration is invariant under algebra automorphisms. In \autoref{section: analysis of stability} we use these properties to characterize the semistable algebras. We also show that to understand the canonical filtration for Lie algebras one can reduce to the nilpotent case, and that the canonical filtration for Lie algebras is a filtration by ideals. We also compute the canonical filtrations of some algebras constructed from semisimple associative algebras. In \autoref{section: strategy to compute canonical filtrations}, we give a strategy to reduce the set of filtrations one must consider when computing the canonical filtration of $\bZ^r$-graded algebras. When the homogeneous components have dimension one, this reduces the computation to a single optimization problem. In \autoref{section: examples and computations}, we use this strategy to compute the canonical filtration in some explicit examples, such as monomial algebras, algebras of upper triangular matrices, filiform Lie algebras, and Borel subalgebras of semisimple Lie algebras over $\bC$. In \autoref{section: characterizations of graded-semistable algebras}, we obtain results on the algebraic structure of associative graded-semistable algebras. In particular, this gives some results about the structure of the associated graded algebra of the canonical filtration of an associative algebras. We also show that the operator defining the canonical filtration of a graded-semistable algebra is dual to a certain functional on the space of $\bZ$-gradings.

\noindent \textbf{Acknowledgements}: The author would like to thank Daniel Halpern-Leistner for his guidance, suggestions, and helpful discussions. We also thank Ritvik Ramkumar for his helpful feedback on drafts of this paper.

\section{Preliminaries} \label{section: preliminaries}

We fix a perfect field $\basefield$ and let $V$ be a finite-dimensional $\basefield$-vector space of dimension $d$. When no confusion is likely to arise, we will simply refer to an associative or Lie algebra as an algebra.

\subsection{Schemes of Algebra Structures}

We can identify $\Hom(V^{\otimes 2}, V)$ with the affine scheme $\Spec(\Sym^{\bullet}(V^{\otimes 2} \otimes V^{\dual}))$, where $V^{\dual} = \Hom(V, \basefield)$. Choose a basis $ \cB = \{ x_1, \ldots, x_{d} \}$ of $V$. Any bilinear map $\mu \in \Hom(V^{\otimes 2}, V)$ is completely determined by structure constants $c_{ij}^l \in \basefield$, $1 \leq i,j,l \leq d$, as follows. For $x_i, x_j \in \cB$ we can write
\[ \mu(x_i, x_j) = \sum_{l = 1}^d c_{ij}^l x_l \]
and conversely any set of $c^l_{ij}$ as above determines a bilinear map. To give a (not necessarily unital) associative algebra structure on $V$ requires that $\mu$ satisfies the associativity condition $\mu \circ (\mu \otimes \id_V) = \mu \circ (\id_V \otimes \mu)$. Similarly, to give $V$ a Lie algebra structure requires that $\mu$ is alternating and satisfies the Jacobi identity. These conditions are polynomial in the structure constants $c_{ij}^l$, hence the bilinear maps defining an algebra structure along with a choice of basis are a closed subscheme of $\Hom(V^{\otimes 2}, V)$. Thus we can make the following definition.

\begin{defn} \label{defn: schemes of associative and Lie algebras}
The schemes $\ModAlg$ and $\ModLieAlg$ are the closed subschemes of $\Hom(V^{\otimes 2}, V)$ parameterizing associative and Lie algebra structures on $V$ plus the choice of a basis of $V$.
\end{defn}

Observe that two choices of bases and structure constants may define isomorphic algebras. To rectify this, we quotient out by all possible change of bases. More precisely, let $Y$ be either $\ModAlg$ or $\ModLieAlg$. It follows from the universal properties of $\Spec$, $\Sym$, and the tensor product that $Y$ represents the functor which assigns to a $\basefield$-algebra $R$ the set of algebra structures on $V \otimes R$ of the type corresponding to $Y$. Let $R$ be a $\basefield$-algebra and let $\mu$ be an $R$-bilinear map on $V$. Then the set of $R$-linear automorphisms $\GL(V)(R) = \Aut_R(V \otimes R)$ acts on $V \otimes R$ by $g \cdot \mu := \mu_g$, where
\[ \mu_g(x, y) = g \mu(g^{-1} x, g^{-1} y) \]
This is functorial in $R$, and hence defines an action of $\GL(V)$ on $\Hom(V^{\otimes 2}, V)$. A straightforward computation shows that if $\mu$ defines an associative or Lie algebra structure on $V \otimes R$ then so does $\mu_g$. It is also straightforward to show that the linear map $\ell_g \in \End_R(V \otimes R)$, defined as $\ell_g(x) = g x$, is an isomorphism of $\mu$ and $\mu_g$ in each of the cases above. Thus the $\GL(V)$-orbits of each of the subschemes $X$ in definition (\ref{defn: schemes of associative and Lie algebras}) corresponds to isomorphic algebra structures.

\begin{remark}
In the future we will abuse notation and a point $\mu$ in $\ModAlg$ or $\ModLieAlg$ giving an algebra structure on $V$ by a pair $(A, \mu_A)$, $A$, or $\mu_A$ when there is no risk of confusion.
\end{remark}

\subsection{Geometric Invariant Theory on the Parameter Schemes}

In this section, we define a notion of semistability via $\GL(V)$-linearized line bundles on $\ModAlg$ and $\ModLieAlg$; let $Y$ be either of these schemes. We take the trivial line bundle $Y \times \bA^1 \to Y$, and linearize it by lifting the $\GL(V)$ action via a character $\rho = \det^{p}$. In coordinates, this action is given by
\[ g \cdot (\mu_A, z) = (g \cdot \mu_A, \rho(g) z ) \]
where $z$ is the coordinate of the line bundle. The stability of an algebra will depend only on the sign of $p$, so there are only two non-trivial choices to consider: $\det$ and $\det^{-1}$. It will turn out that only one of these yields an interesting stability condition.

The definition of semistability can be stated in terms of certain invariant sections of the linearized line bundle, but in practice this is difficult to work with. For our definition we will use an equivalent criterion, which was originally due to Mumford \cite{mumford1994geometric} and adapted to the case of affine schemes by others (see for instance \cite{hoskins2014stratifications}, \cite{king1994moduli}). Before stating this criterion, we recall the definition of the pairing between one-parameter subgroups and characters for the special cases of $\GL(V)$ and a split torus $T \subseteq \GL(V)$. Given a split torus $T$ of rank $d$, we have the following identifications of the group of one-parameter subgroups and the group of characters of $T$, respectively:
\[ \Gamma(T) = \Hom(\bG_m, T) \cong \bZ^d, \quad X(T) = \Hom(T, \bG_m) \cong \bZ^d \]
Given a one-parameter subgroup $\lambda$ of $T$ and a character $\rho$, the composition $\rho \circ \lambda : \bG_m \to \bG_m$ is classified by the integer $(\rho, \lambda)$ such that $(\rho \circ \lambda)(t) = t^{(\rho, \lambda)}$. Under the previous identifications, this integer is just the usual dot product on $\bZ^d$: if $\lambda$ and $\rho$ are identified with vectors $(m_1, \ldots, m_d)$ and $(\rho_1, \ldots, \rho_d)$ in $\bZ^d$ respectively, then
\[ (\rho, \lambda) = \sum_{i = 1}^d \rho_i m_i \]
If $\lambda : \bG_m \to \GL(V)$ is a one-parameter subgroup such that the image of $\lambda(t)$ is contained in some split maximal torus $T$, we can restrict the character $\rho = \det^p$ to $\bG_m$ along $\lambda$. Using the identifications above, we have
\[ (\rho, \lambda) = p \sum_{i=1}^d m_i \]
for integers $m_i \in \bZ^d$ corresponding to $\lambda$. With this in place we can state a criterion for testing semistability.

\begin{thm}[Affine Hilbert-Mumford Criterion] \label{thm: Affine Hilbert-Mumford criterion}
Let $Y$ be one of $\ModAlg$ or $\ModLieAlg$ and let $\rho$ be a character of $\GL(V)$. A point $\mu_A \in Y$ is semistable if and only if for all one-parameter subgroups $\lambda : \bG_m \to \GL(V)$ such that $\lim_{t \to 0} \lambda(t) \cdot \mu_A$ exists in $Y$, we have $(\rho, \lambda) \geq 0$. Otherwise we say that $\mu_A$ is unstable, and if $(\rho, \lambda) < 0$ for some $\lambda$ we say that $\lambda$ is destabilizing for $\mu_A$.
\end{thm}

\subsection{Algebra Filtrations} \label{section: algebra filtrations}

In this section we define our notion of an algebra filtration. This will enable us to interpret the existence of the limit in the Hilbert-Mumford criterion (\ref{thm: Affine Hilbert-Mumford criterion}) explicitly and eventually define the canonical filtration. We will also give two other equivalent ways to define such a filtration, which will be useful for later computations. We will also recall the associated graded algebra of a filtration.

\begin{defn} \label{defn: algebra filtration}
An algebra filtration $F_\bullet$ of an algebra $A$ is a sequence of $\basefield$-subspaces $F_m$, $m \in \bZ$, such that
\begin{enumerate}
    \item[(i)] $F_i \supseteq F_j$ when $i \leq j$.
    \item[(ii)] $\mu_A(F_i, F_j) \subseteq F_{i + j}$ for all $i,j$.
    \item[(iii)] $F_i = A$ for $i \ll 0$ and $F_i = 0$ for $i \gg 0$.
\end{enumerate}
The trivial filtration is the filtration $F_\bullet$ such that $F_n = A$ for all $n \leq 0$ and $F_n = 0$ for $n > 0$.
\end{defn}

When no confusion is likely to arise, we will simply refer to an algebra filtration of an algebra $A$ as a filtration. It turns out that the limit in (\ref{thm: Affine Hilbert-Mumford criterion}) will exist precisely when the one-parameter subgroup $\lambda(t)$ defines an algebra filtration.

\begin{prop} \label{prop: one-parameter subgroups are filtrations of algebras}
Let $A$ be an associative (resp. Lie) algebra represented by a point $\mu_A$ in $\ModAlg$ (resp. $\ModLieAlg$). To give a one-parameter subgroup $\lambda(t) : \bG_m \to \GL(V)$ such that $\lim_{t \to 0} \lambda(t) \cdot \mu_A$ exists is the same as giving an algebra filtration $F(\lambda)_\bullet$ of $A$.
\end{prop}

\begin{proof}
Since $\bG_m$ is linearly reductive, there is a basis $e_1, \ldots, e_{d}$ for which the induced $\bG_m$-action is diagonal and
\[ A = \bigoplus_{m \in \bZ} A_m \]
where $A_m = \Span\{ a \in A \mid \lambda(t) \cdot a = t^{m} a \}$. We obtain a vector space filtration of $A$, which evidently satisfies (i) and (iii), by setting $F(\lambda)_n = \bigoplus_{n \geq 0} A_n$. Consider the underlying bilinear map $\mu_A$ of $A$, written as a tensor in this diagonalizing basis as
\[ \mu_A = \sum_{i,j,k} c_{ij}^k \cdot e_i^\vee \otimes e_j^\vee \otimes e_k \]
where $e_i^\vee$ denotes the vector in $V^\vee$ dual to $e_i$. Observe that the induced action of $\bG_m$ on $\Hom(V^{\otimes 2}, V)$ is given by 
\[ \lambda(t) \cdot \mu_A = \sum_{i,j,k} t^{m_k - m_i - m_j} c_{ij}^k \cdot e_i^\vee \otimes e_j^\vee \otimes e_k \]
The limit as $t \to 0$ exists if and only if $m_k - m_i - m_j \geq 0$ for all $i, j, k$, which implies that $\mu_A(F(\lambda)_{m_i}, F(\lambda)_{m_j}) \subseteq F(\lambda)_{m_k}$. This gives (ii). Conversely, given $F_\bullet$ satisfying (i) - (iii) we can choose a basis as above and define a one-parameter subgroup by declaring that $\lambda(t) e_i = t^{m_i} e_i$ for the unique $m_i \in \bZ$ such that $e_i \in F_{m_i} \setminus F_{m_i + 1}$. The condition $\mu_A(F_i, F_j) \subseteq F_{i+j}$ then implies that the limit as $t \to 0$ exists.
\end{proof}

If we wish to emphasize the one-parameter subgroup $\lambda(t)$ from which a filtration arises, we shall write $F(\lambda)_\bullet$, but otherwise we will simply write $F_\bullet$. Next we give an equivalent characterization of filtrations in terms of weight functions.

\begin{defn} \label{defn: weight function}
For an algebra $A$, a function $w : A \to \bR \cup \{ \infty \}$ is called a weight function if for all $a, b \in A$ and $c \neq 0 \in \basefield$,
\begin{enumerate}
    \item[(i)] $w(c a) = w(a)$.
    \item[(ii)] $w(\mu_A(a,b)) \geq w(a) + w(b)$.
    \item[(iii)] $w(a + b) \geq \min(w(a), w(b))$
    \item[(iv)] $w(a) = \infty \iff a = 0$
\end{enumerate}
\end{defn}

\begin{prop} \label{prop: filtrations = weight functions}
To give an algebra filtration $F_\bullet$ of $A$ is the same as giving a weight function on $A$.
\end{prop}

\begin{proof}
Given a weight function $w$, we define a filtration 
$F(w)_\bullet$ where 
\[ F(w)_m = \{ a \in A \mid w(a) \geq m \} \]
Conditions (i), (iii), and (iv) of weight functions imply that $F(w)_m$ is a subspace of $A$. Condition (ii) shows that $\mu( F(w)_m, F(w)_n ) \subseteq F(w)_{m + n}$. If we choose a basis $e_1, \ldots, e_d$ for $A$, then for any $x \neq 0 \in A$ we have
\[ w(x) = w \left( \sum_{i = 1}^d c_i e_i \right) \geq \min_{c_i \neq 0}(w(e_i)) \]
from condition (iii), which shows that there is some $N$ such that $w(x) \geq N$ for any $x \neq 0$. Therefore, $F(w)_n = A$ for $n \ll 0$. Since $A$ is finite-dimensional and the $F(w)_\bullet$ defines a decreasing sequence of subspaces, there is some $N \gg 0$ such that $F(w)_n = F(w)_N$ for $n \geq N$. But $x \in F(w)_N$ implies that $x \in F(w)_n$ for all $n$, hence $x = 0$ by condition (iv). Thus, $F(w)_n = 0$ for $n \gg 0$, which shows that $F(w)_\bullet$ is bounded. Hence $F(w)_\bullet$ is an algebra filtration.

Conversely, given an algebra filtration $F_\bullet$ we obtain a weight function $w_{F_\bullet}$ by setting 
\[ w_{F_\bullet}(a) = \max( \{ n \in \bZ \mid a \in F_n \} ) \]
if the maximum exists and $w_{F_\bullet}(a) = \infty$ otherwise. Condition (i) of the weight function follows from the fact that each $F_n$ is a subspace. Since $w_{F_\bullet}(a) = \infty$ if and only if $a \in F_n$ for all $n$ if and only if $a = 0$, we see that condition (iv) of the weight function holds. If we set $n = w_{F_\bullet}(a)$ and $m = w_{F_\bullet}(b)$, then we have $a \in F_n$ and $b \in F_m$, hence $\mu_A(a, b) \in F_{n +m}$. It follows from the definition that 
\[ w_{F_\bullet}(\mu_A(a,b)) \geq n + m = w_{F_\bullet}(a) + w_{F_\bullet}(b) \]
which gives condition (ii). Now if $a, b \in A$ and $a + b = 0$ then (iii) holds trivially, so assume that $a + b \neq 0$. Without loss of generality, we can assume that $n := w_{F_\bullet}(a) \leq w_{F_\bullet}(b)$. Then we have $a, b \in F_n$, hence $a + b \in F_n$, and it follows that
\[ w_{F_\bullet}(a + b) \geq n = \min(w_{F_\bullet}(a), w_{F_\bullet}(a)) \]
This shows that condition (iii) holds, so $w_{F\bullet}$ is a weight function.
\end{proof}

In later computations we will want to start with a $\basefield$-basis $\cB$ of an algebra $A$ and assign weights to for each $x \in \cB$ in such a way that there is a weight function $w$ which agrees with this assignment, and hence determines a filtration. We define the following key data for this construction.

\begin{defn} \label{defn: weighted algebra basis}
A basis $\cB = \{ a_1, \ldots, a_d \}$ with weights $w_1, \ldots, w_d \in \bR$ of an algebra $A$ is said to be a weighted algebra basis if it satisfies the weight inequality conditions: for $x,y \in \cB$, if
\[ 
\mu(x,y) = \sum_{a \in \cB} c_a \cdot a, \quad c_a \in \basefield 
\]
then
\begin{equation} \label{eqn: weight basis inequality condition}
w(a) < w(x) + w(y) \implies c_a = 0
\end{equation}
\end{defn}

Given a weighted algebra basis $\cB$ of $A$ with weights $\{ w(a) \mid a \in \cB \}$, one can define a weight function by finding the basis expansion of $x \in A$ 
\[ x = \sum_{a \in \cB} c_a \cdot a \]
and then setting $w(x) = \min \{ w(a) \mid c_a \neq 0 \}$. 

The construction above allows us to construct filtrations of an algebra $A$ from a basis, but it will also be useful to build algebra filtrations of $A$ from the data of a vector space filtration of $A$. This brings us to our final equivalent characterization algebra filtrations. We first make the following preliminary definition. 

\begin{defn} \label{defn: weighted vector space filtration}
Let $V$ be a vector space over $\basefield$. A $\bZ$-weighted vector space filtration of $V$ is the data of finite filtration of subspaces
\[ V = G_{(0)} \supsetneq G_{(1)} \supsetneq \ldots \supsetneq G_{(n)} \supsetneq G_{(n+1)} = 0 \]
and a weight $w_i \in \bZ$ for each $G_i$ with $0 \leq i \leq n$, and $w_{n+1} = \infty$. We denote this data by $(G_{(\bullet)}, w_\bullet)$.
\end{defn}

An algebra filtration is then just a weighted vector space filtration of $A$ subject to certain constraints on the weights.

\begin{prop} \label{prop: algebra filtration from weighted filtration}
To give a filtration $F_\bullet$ of an algebra $A$ is equivalent to giving a weighted vector space filtration $(G_{(\bullet)}, w_\bullet)$ satisfying the weight condition: if $p$ is the minimal index such that $w_i + w_j \leq w_p$, then $\mu_A(G_{(i)}, G_{(j)}) \subseteq G_{(p)}$.
\end{prop}

\begin{proof}
To avoid cases and make the proof more uniform, we set $F_{\infty} = 0$ and use the convention that $\infty + m = m + \infty = \infty$ for all $m \in \bZ$. Starting with the data of an algebra filtration $F_\bullet$, we set
\[ \{ w_0, \ldots, w_n \} := \{ m \in \bZ \mid F_m \neq F_{m+1} \} \]
with $w_i < w_{i + 1}$, and set $w_{n+1} = \infty$. Note that the number of $w_i$ is finite since $F_\bullet$ is a bounded filtration. Set $G_{(i)} = F_{w_i}$. We will show that the weight condition is satisfied; assume $w_p \geq w_i + w_j$ with $p$ minimal. We have
\[ \mu_A(G_{(i)}, G_{(j)}) = \mu_A(F_{w_i}, F_{w_j}) \subseteq F_{w_i + w_j} \]
Observe that if $r \in \bZ$ is such that $w_{p-1} < r \leq w_p$, then $F_r = F_{w_p}$. From this observation and the minimality of $p$ it follows that $F_{w_i + w_j} = F_{w_p} = G_{(p)}$, as required.

Conversely, suppose we have a weighted vector space filtration satisfying the weight condition. For each $i \in \bZ$, denote by $j(i)$ the minimal index such that $i \leq w_{j(i)}$ and set $F_i = G_{(j(i))}$. We show that conditions (i) - (iii) of definition (\ref{defn: algebra filtration}) are satisfied. Note that by construction we have $F_i = A$ when $j(i) \leq w_0$ and $F_i = 0$ if $i > w_n$, so condition (iii) is satisfied. Since $i < i + 1 \leq w_{j(i+1)}$ by the definition of $j(i+1)$, we have $j(i) \leq j(i+1)$ from the minimality of $j(i)$. Therefore,
\[ F_i = G_{(j(i))} \supseteq G_{(j(i + 1))} = F_{i+1} \]
and hence condition (i) is satisfied. For $s, t \in \bZ$, we have 
\[ F_s = G_{(j(s))}, \quad F_t = G_{(j(t))}, \quad F_{s + t} = G_{(j(s + t))} \]
Choose the minimal $p$ such that $w_{j(s)} + w_{j(t)} \leq w_p$. Then we have
\[ s + t \leq w_{j(s)} + w_{j(t)} \leq w_p \]
so $w_{j(s + t)} \leq w_p$ by the minimality of $j(s + t)$, hence $G_{(p)} \subseteq G_{(j(s + t))}$. Therefore,
\[ \mu_A(F_s, F_t) = \mu_A(G_{(j(s))}, G_{(j(t))}) \subseteq G_{(p)} \subseteq G_{(j(s + t))} = F_{s + t} \]
which shows that condition (ii) holds and completes the proof.
\end{proof}

Lastly, we recall the definition of the associated graded algebra of a filtration:

\begin{defn} \label{defn: associated graded of a filtration}
Let $F_\bullet$ be a filtration of an algebra. The associated graded algebra of $F_\bullet$ is
\[ \Gr(F_\bullet) = \Gr(F_\bullet) := \bigoplus_{m} F_m / F_{m + 1} \]
The multiplication for $\Gr(F_\bullet)$ is given on homogeneous elements $u \in F_m / F_{m + 1}$ and $v \in F_n / F_{n + 1}$ by choosing lifts of $x$ and $y$ in $F_m$ and $F_n$ of $u$ and $v$, respectively, and then defining $uv$ to be the class of $xy \in F_{m + n} / F_{m + n + 1}$. 
\end{defn}

\subsection{The Hilbert-Mumford Weight and the Numerical Invariant}

We can now explicitly interpret the pairing between one-parameter subgroups and characters appearing in the Hilbert-Mumford criterion thanks to (\ref{prop: one-parameter subgroups are filtrations of algebras}). Given a one-parameter subgroup $\lambda(t)$, we can choose a basis in which the induced action of $\bG_m$ is diagonal so that in these coordinates $\lambda(t) = \diag(t^{m_1}, \ldots, t^{m_d})$. If $F_\bullet = F(\lambda)_\bullet$ is the filtration corresponding to $\lambda$, then $\dim(F_m / F_{m+1})$ is the number of times $m$ appears as one of the $m_i$, and thus for the character $\rho = \det^p$ we have
\[ (\rho, \lambda) = p \sum_{i=1}^d m_i = p \sum_{m \in \bZ} m \dim(F_m / F_{m+1}) \] 
Note that the right-hand sum is finite as only finitely many of the quotients $F_m / F_{m+1}$ are non-zero. It is now clear that the sign of $(\rho, \lambda)$ depends only on the sign of $p$, hence we need only consider the characters $\det^{\pm 1}$. If we consider the filtration $F_\bullet$ given by $F_n = A$ if $n < 0$ and $F_n = 0$ if $n \geq 0$ and its corresponding one-parameter subgroup $\lambda$, then for $\rho = \det$ we have
\[ (\rho, \lambda) = \sum_{m \in \bZ} m \dim(F_m / F_{m+1}) = - \dim(A) < 0 \]
so \emph{every} algebra is unstable with respect to this character. Hence, we will only consider the character $\det^{-1}$. See also remark (\ref{remark: stratifications for the det character}) for the effect this choice has on the canonical filtrations.

If an algebra $A$ is unstable, a theorem of Kempf \cite{kempf1978instability} shows that there is a class of optimal destabilizing one-parameter subgroups which destabilizes $A$ ``most rapidly''. To define the optimality condition, we also require a norm on one-parameter subgroups. Since one-parameter subgroups correspond to filtrations in our situation, we will define this norm in terms of filtrations as follows.

\begin{defn} \label{defn: weight and norm of a filtration}
Given an algebra $A$ and a filtration $F_\bullet$, we set
\[  \wt(A, F_\bullet) = \sum_{m \in \bZ} m \dim(F_m / F_{m+1}), \quad || F_\bullet || = \sum_{m \in \bZ} m^2 \dim(F_m / F_{m+1}) \]
\end{defn}

Note that the norm in Kempf's theorem is required to be invariant under the Weyl group of a fixed maximal torus $T$. This is indeed the case for the norm above: one can check that the norm above corresponds to the usual norm on one-parameter subgroups of $T$ induced by the identification of one-parameter subgroups and $\bZ^r$, which has the required property. With all this in place, we can finally define our notion of semistability:

\begin{defn} \label{defn: semistablity for algebras}
We say an algebra $A$ is semistable (with respect to the $\det^{-1}$ character) if for all filtrations $F_\bullet$ we have $ - \wt(A, F_\bullet) \geq 0$.
\end{defn}

Next we define a numerical invariant on our parameter schemes, which is the function which determines the class of optimal destabilizing one-parameter subgroups for a given algebra.

\begin{defn} \label{defn: numerical invariants}
The numerical invariant is
\[ \nu(A, F_\bullet) = \frac{- \wt(A, F_\bullet)}{||F^\bullet||} \]
where $A$ is an (associative or Lie) algebra and $F_\bullet$ is a filtration of $A$.
\end{defn}

For a given algebra, the optimal class of one-parameter subgroups all define the same filtration (see \cite[Theorem 3.4]{kempf1978instability}), so we can make the following definition.

\begin{defn} \label{defn: the canonical filtration}
The canonical filtration $F(A)^{\Can}_\bullet$ of an algebra $A$ is the filtration which minimizes the numerical invariant $\nu$ if $A$ is unstable and is the trivial filtration if $A$ is semistable.
\end{defn}

\begin{remark} \label{remark: stratifications for the det character}
Although the notion of semistability arising from the $\det$ character is trivial (all algebras are unstable), one can ask if the canonical filtrations determined by the corresponding numerical invariant $\frac{\wt(A, F_\bullet)}{||F_\bullet||}$ are interesting. It turns out that this is not the case: one can show that for this alternate numerical invariant, the canonical filtration of any algebra $A$ is given by the filtration $F_\bullet$ where $F_n = A$ if $n < 0$ and $F_n = 0$ if $n \geq 0$.
\end{remark}

\section{Properties of the Canonical Filtration} \label{section: properties of the canonical filtration}

Given filtrations $F(A_1)_\bullet, F(A_2)_\bullet$ of algebras $A_1, A_2$, we define the direct sum filtration $F(A_1)_\bullet \oplus F(A_2)_\bullet$ of $A_1 \oplus A_2$ to be the filtration such that the $m$-th filtered piece is $F(A_1)_m \oplus F(A_2)_m$. For a filtration $G_\bullet$ of an algebra and $c \in \bZ_{> 0}$, we denote by $G_{c \bullet}$ the filtration $H_\bullet$ such that $H_n = G_{c n}$.

\begin{remark}
The reader may object to speaking of the direct $A \oplus B$ sum of two associative algebras $A, B$, since if $A$ and $B$ are unital then the canonical inclusion maps of vector spaces are not unital algebra homomorphisms. However, we have decided to use direct sum notation since we do not insist that associative algebras are unital, our arguments do not use unital homomorphisms in a critical way, and the use of direct sums is standard for Lie algebras.
\end{remark}

\begin{prop} \label{prop: canonical filtration for direct sum of algebras}
Let $A = A_1 \oplus A_2$ be a direct sum of algebras and set
\[ \ell_i = \wt(A_i, F(A_i)^{\Can}_\bullet), \quad b_i = ||F(A_i)^{\Can}_\bullet||^2 \]
where $F(A_i)^{\Can}_\bullet$ is the canonical filtration of $A_i$. 
\begin{enumerate}
    \item[(i)] If both $A_i$ are semistable then so is $A$, and hence the canonical filtration of $A$ is the trivial filtration.
    \item[(ii)] If $A_1$ is semistable and $A_2$ is unstable then
    \[ F(A)^{\Can}_\bullet = F(A_1)^{\Can}_\bullet \oplus F(A_2)^{\Can}_{\bullet} \]
    \item[(iii)] If both $A_1$ and $A_2$ are unstable then
    \[ F(A)^{\Can}_\bullet = F(A_1)^{\Can}_{b_2 \ell_1 \bullet} \oplus F(A_2)^{\Can}_{b_1 \ell_2 \bullet} \]
\end{enumerate}
\end{prop}

\begin{proof}
We set
\[ F^1_\bullet = F(A)^{\Can}_\bullet \cap A_1, \quad F^2_\bullet = (F(A)^{\Can}_\bullet + A_1) / A_1 \]
and note that a routine check shows that $F^i_\bullet$ is an algebra fitlration of $A_i$. Observe that since
\[ (F(A)^{\Can}_m + A_1) / A_1 \cong F(A)^{\Can}_m / (A_1 \cap F(A)^{\Can}_m) \]
as vector spaces, the additivity of dimension implies that
\begin{equation} \label{eqn: filtration weight sum}
\wt(A, F(A)^{\Can}_\bullet) = \wt(A_1, F^1_\bullet) + \wt(A_2, F^2_\bullet) 
\end{equation}
\begin{equation} \label{eqn: filtration norm sum}
||F(A)^{\Can}_\bullet|| = \sqrt{||F^1_\bullet||^2 + ||F^2_\bullet||^2}
\end{equation}
We handle each of the cases in the proposition statement separately.

(i): From equation (\ref{eqn: filtration weight sum}) and the assumption that both $A_1$ and $A_2$ are semistable, we have
\[ - \wt(A, F(A)^{\Can}_\bullet) = - [ \wt(A_1, F^1_\bullet) + \wt(A_2, F^2_\bullet) ] \geq 0 \]
This is only possible if $A$ is semistable, which implies the claim.

(ii): Set 
\[ G_\bullet = F(A_1)^{\Can}_\bullet \oplus F(A_2)^{\Can}_\bullet \]
Since $F(A_1)^{\Can}_\bullet$ is the trivial filtration in this case, it is immediate that
\[ \numinvariant(A, G_\bullet) = \numinvariant(A_2, F(A_2)^{\Can}_{\bullet}) \]
We also have that
\begin{equation} \label{eqn: direct sum inequality}
\frac{- \wt(A, G_\bullet)}{||G_\bullet||} \geq \frac{- \wt(A, F(A)^{\Can}_\bullet)}{||F(A)^{\Can}_\bullet||}
\end{equation}
The assumption that $A_1$ is semistable as well as the equalities in equations (\ref{eqn: filtration weight sum}), (\ref{eqn: filtration norm sum}) implies that
\[
\frac{- \wt(A, F(A)^{\Can}_\bullet)}{|| F(A)^{\Can}_\bullet||} \geq \frac{- \wt(A_2, F^2_\bullet)}{||F^2_\bullet||} 
\geq \frac{- \wt(A_2, F(A_2)^{\Can}_{ \bullet})}{||F(A_2)^{\Can}_{\bullet}||} = \frac{- \wt(A, G_\bullet)}{|| G_\bullet||} \]
This shows that equality holds in (\ref{eqn: direct sum inequality}). Taking square roots and negating shows that $G_\bullet$ is the canonical filtration.

(iii): First note that $\ell_i$ and $b_i$ are both positive. 
Set 
\[ G_\bullet = F(A_1)^{\Can}_{b_2 \ell_1 \bullet} \oplus F(A_2)^{\Can}_{b_1 \ell_2 \bullet} \]
By Titu's lemma, the equalities in equations (\ref{eqn: filtration weight sum}) and (\ref{eqn: filtration norm sum}), and the definition of the canonical filtrations, we find that
\begin{align*}
\frac{\wt(A, F(A)^{\Can}_\bullet)^2}{||F(A)^{\Can}_\bullet||^2} &\leq \frac{\wt(A_1, F^1_\bullet)^2}{||F^1_\bullet||^2} + \frac{\wt(A_2, F^2_\bullet)^2}{||F^2_\bullet||^2} \\
&\leq \frac{\wt(A_1, F^{\Can}(A_1)_\bullet)^2}{||F^{\Can}(A_1)_\bullet||^2} + \frac{\wt(A_2, F^{\Can}(A_2)_\bullet)^2}{||F^{\Can}(A_2)_\bullet||^2} \\
&= \frac{\wt(A_1, F^{\Can}(A_1)_{b_2 \ell_1 \bullet})^2}{||F^{\Can}(A_1)_{b_2 \ell_1 \bullet}||^2} + \frac{\wt(A_2, F^{\Can}(A_2)_{b_1 \ell_2 \bullet})^2}{|| F^{\Can}(A_2)_{b_1 \ell_2 \bullet}||^2}
\end{align*}
where the last equality follows from the fact that the summands are invariant under scaling the weights of the filtration by the positive numbers $b_2 \ell_1$, $b_1 \ell_2$. Observe that if $H_\bullet$ is a filtration of an algebra $R$ and $c \in \bZ_{\geq 0}$, then $\wt(R, H_{c \bullet}) = c \wt(R, H_{\bullet})$ and $|| H_{c \bullet}|| = c ||H_{\bullet}||$. Thus
\[ \frac{\wt(A_1, F^{\Can}(A_1)_{b_2 \ell_1 \bullet})}{||F^{\Can}(A_1)_{b_2 \ell_1 \bullet}||^2} = \frac{(b_2 \ell_1) \ell_1}{(b_2 \ell_1)^2 b_1} = \frac{1}{b_1 b_2} \]
\[ \frac{\wt(A_2, F^{\Can}(A_2)_{b_1 \ell_2 \bullet})}{||F^{\Can}(A_2)_{b_1 \ell_2 \bullet}||^2} = \frac{(b_1 \ell_2) \ell_2}{(b_1 \ell_2)^2 b_2} = \frac{1}{b_1 b_2} \]
This implies equality must hold in Titu's lemma, so we have
\[ \frac{\wt(A, G_\bullet)^2}{||G_\bullet||^2} = \frac{\wt(A_1, F^{\Can}(A_1)_{b_2 \ell_1 \bullet})^2}{|| F^{\Can}(A_1)_{b_2 \ell_1 \bullet}||^2} + \frac{\wt(A_2, F^{\Can}(A_2)_{b_1 \ell_2 \bullet})^2}{||F^{\Can}(A_2)_{b_1 \ell_2 \bullet}||^2} \]
Combining all of the above, we find that
\[ \frac{\wt(A, G_\bullet)^2}{||G_\bullet||^2} \leq \frac{\wt(A, F(A)^{\Can}_\bullet)^2}{||F(A)^{\Can}_\bullet||^2} \leq \frac{\wt(A, G_\bullet)^2}{||G_\bullet||^2} \]
and equality holds. Taking square roots and negatives, we see that $G_\bullet$ is the canonical filtration.
\end{proof}

We now come to some of the most important properties of canonical filtrations. In subsequent sections, we shall use the following proposition to determine the semistable locus of $\ModAlg$ and $\ModLieAlg$ and reduce the computation of canonical filtrations to tractable convex optimization problems.

\begin{prop} \label{prop: action of automorphisms and derivations on canonical filtrations}
For any algebra $A$ and any $n \in \bZ$, $F(A)^{\Can}_n$ is stable under the action of $\Aut(A)$ and $\Der(A)$.
\end{prop}

\begin{proof}
Note that the statement is trivial when $A$ is semistable, so we can assume that $A$ is unstable. For any $\varphi \in \Aut(A)$, $\varphi(F(A)^{\Can}_\bullet)$ is also a filtration of $A$ and $\dim \varphi(F(A)^{\Can}_n) = \dim F(A)^{\Can}_n$, which implies that $\numinvariant(A, \varphi(F(A)^{\Can}_\bullet)) = \numinvariant(A, F(A)^{\Can}_\bullet)$. By uniqueness, $\varphi(F(A)^{\Can}_\bullet)$ is the canonical filtration up to some positive scaling of the weights. The set of indexes $n$ where $\varphi(F(A)^{\Can}_n) \neq \varphi(F(A)^{\Can}_{n+1})$ is the same as the set of $n$ where $F(A)^{\Can}_n \neq F(A)^{\Can}_{n+1}$ since $\varphi$ preserves dimension, so the scaling factor must be $1$. This shows that $\varphi(F(A)^{\Can}_n) = F(A)^{\Can}_n$. Since $\Aut(A)$ is an algebraic group (\cite[Exercise 7.3]{humphreys2012linear}) and $\Der(A)$ is the Lie algebra of $\Aut(A)$, (\cite[Exercise 10.1]{humphreys2012linear}) shows that if a subspace $W$ of $A$ is stable under $\Aut(A)$ then $W$ is also stable under $\Der(A)$.
\end{proof}

\begin{lemma} \label{lemma: rationality of the canonical filtration}
Let $\overline{\basefield}$ denote the algebraic closure of $\basefield$, and denote by $A_{\overline{\basefield}}$ the base change of the $k$-algebra $A$ by $\overline{\basefield}$. Then $F(A)^{\Can}_\bullet \otimes_{\basefield} \overline{\basefield} = F(A_{\overline{\basefield}})^{\Can}_\bullet$. As a consequence, if for some $n$ we have $x \otimes 1 \in F(A_{\overline{\basefield}})^{\Can}_n$, then $x \in F(A)^{\Can}_n$.
\end{lemma}

\begin{proof}
The first claim follows from \cite[Theorem 4.2]{kempf1978instability}, which in this context says that the canonical filtration of $A_{\overline{\basefield}}$ is defined over $k$. For the second claim, let $W = \Span_{\basefield}(x)$ and consider the inclusion $\iota: W \cap F(A)^{\Can}_n \to W$. After base-changing to $\overline{\basefield}$ and using the fact that tensor products commute with finite intersections, our assumptions imply that $\iota \otimes \overline{\basefield}$ is an isomorphism. Since base change by $\overline{\basefield}$ is faithfully flat, $\iota$ must be an isomorphism and hence $x \in F(A)^{\Can}_n$.
\end{proof}

\begin{coroll} \label{coroll: canonical filtration is a direct sum of graded pieces}
Suppose that there is a morphism of algebraic groups $T \to \Aut_{\basefield}(A)$, where $T$ is a split torus. Then for every $m$,
\[ F(A)^{\Can}_m = \bigoplus_{\chi} (A_{\chi} \cap F(A)^{\Can}_m)  \]
where the sum is over all characters $\chi$ of $T$ and $A_{\chi} := \{ a \mid t \cdot a = \chi(t) a \}$.
\end{coroll}

\begin{proof}
After base changing and applying lemma (\ref{lemma: rationality of the canonical filtration}), we may assume that $\basefield$ is algebraically closed and hence infinite. Since $T$ is linearly reductive, for $a \neq 0 \in F(A)^{\Can}_m$ we may write $a = \sum_{i=1}^{s} a_i$, where $a_i \neq 0 \in a_i$ and the $\chi_i$ are distinct characters. The claim is equivalent to showing that $a_i \in F(A)^{\Can}_m$ for all $i$, which we show by induction on $s$. If $s = 1$ there is nothing to show, so assume that $s > 1$. Since the characters $\chi_i$ are distinct, we can choose some $t \in \basefield$ so that $\chi_1(t) \neq \chi_2(t)$. Then we have
\[ t\cdot a - \chi_1(t) a = \sum_{i=2}^s (\chi_i(t) - \chi_1(t)) a_i \in F(A)^{\Can}_m \]
Let $S$ be the set of $i$ such that $\chi_i(t) = \chi_1(t)$. By assumption, $2 \notin S$. Furthermore, this implies that $t \cdot a - \chi_1(t) a \neq 0$, so the induction hypothesis implies that $(\chi_i(t) - \chi_1(t)) a_i \in F(A)^{\Can}_m$ for $i > 1$. If $i \notin S$ then $a_{i} \in F(A)^{\Can}_m$, so $\sum_{i \notin S} a_i \in F(A)^{\Can}_m$. This implies that
\[ \sum_{i \in S} a_{i} = a - \sum_{i \notin S} a_{i} \in F(A)^{\Can}_m \]
If $\sum_{i \in S} a_{i} = 0 $ we are done. Otherwise, since $|S| < s$, another application of the induction hypothesis implies that $a_i \in F(A)^{\Can}_m$ for $i \in S$, concluding the proof.
\end{proof}

Note that to give a $\bZ^r$-grading on $A$ is the same as giving a morphism $T \to \Aut_{\basefield}(A)$ of a split torus $T$. Thus, corollary (\ref{coroll: canonical filtration is a direct sum of graded pieces}) says that for a $\bZ^r$-graded algebra $A$, the filtered piece $F(A)^{\Can}_m$ contains all homogeneous components of its elements.

\begin{coroll} \label{coroll: canonical filtration preserved under simultaneous diagonalization}
Suppose that $\mathfrak{t} \subseteq \Der(A)$ is an Abelian Lie subalgebra consisting of diagonalizable derivations and that $\lambda_i : \mathfrak{t} \to \basefield$ are distinct linear functionals such that
\[ A = \bigoplus_{i = 1}^r E_{\lambda_i} \]
as vector spaces, where $E_{\lambda_i} = \{ a \in A \mid D(a) = \lambda_i(D) a \text{ for all } D \in \mathfrak{t} \}$. Then for all $m$,
\[ F(A)^{\Can}_m = \bigoplus_{i = 1}^r ( F(A)^{\Can}_m \cap E_{\lambda_i} ) \]
\end{coroll}

\begin{proof}
The proof is essentially identical to that of corollary (\ref{coroll: canonical filtration is a direct sum of graded pieces}), using that commuting operators are simultaneously diagonalizable and that derivations preserve the canonical filtration.
\end{proof}

A priori, a filtration $F(A)_\bullet^{\Can}$ is only a filtration by vector spaces. However, when each $F(A)^{\Can}_n$ is an ideal we can say a little more.

\begin{prop} \label{prop: filtration by ideals}
Let $I_\bullet$ be a filtration of an algebra $A$ by ideals. Then $I'_\bullet$, where $I'_n = I_n$ if $n > 0$ and $I'_n = A$ if $n \leq 0$, is also a filtration of $A$. Furthermore, if $\nu(A, I_\bullet) < 0$ then 
\[ \numinvariant(A, I_\bullet) \geq \numinvariant(A, I'_\bullet) \]
Consequently, if the canonical filtration $F(A)^{\Can}_\bullet$ is such that each $F(A)^{\Can}_\bullet$ is an ideal, then $F(A)^{\Can}_n = A$ for $n \leq 0$. In particular, if $A$ is a unital associative algebra then $1 \in F(A)^{\Can}_0 \setminus F(A)^{\Can}_1$ i.e. the unit has weight $0$.
\end{prop}

\begin{proof}
It is immediate from the fact that each $I_n$ is an ideal that $I'_\bullet$ is a filtration. It follows from the construction of $I'_\bullet$ that
\[ \wt(A, I_\bullet) \leq \wt(A, I'_\bullet), \quad ||I_\bullet|| \geq ||I'_\bullet|| \]
The statement about the canonical filtration follows from these inequalities. If $A$ is a unital associative algebra, then since $F(A)^{\Can}_0 = A$ and $1 \not in F(A)^{\Can}_1$ (since this would imply that $1$ is nilpotent), the unit $1$ has weight $0$.
\end{proof}

\begin{coroll} \label{coroll: The canonical filtration of a Lie algebra is an ideal filtration}
Let $A$ be a Lie algebra with canonical filtration $F(A)^{\Can}_\bullet$. Then each subspace $F(A)^{\Can}_n$ is an ideal and $F(A)^{\Can}_n = A$ if $n \leq 0$.
\end{coroll}

\begin{proof}
By proposition (\ref{prop: action of automorphisms and derivations on canonical filtrations}), $F(A)^{\Can}_n$ is stable under the derivation $\ad_x$ for any $x \in A$, so each $F(A)^{\Can}_n$ is an ideal. The second claim now follows from proposition (\ref{prop: filtration by ideals}).
\end{proof}

It is unclear whether the canonical filtration of an associative algebra $A$ is a filtration by ideals. However, the following proposition shows that it suffices to check the weaker condition that for all $n$, $F(A)_{n}^{\Can}$ is closed under multiplication by units.

\begin{prop}
Let $A$ be an associative algebra such that for every unit $u$ and every $n \in \mathbb{Z}$ we have $u F(A)^{\Can}_n, F(A)^{\Can}_n u \subseteq F_\bullet$. Then each $F(A)^{\Can}_n$ is an ideal.
\end{prop}

\begin{proof}
Note that if $x$ is nilpotent then $u = 1 + x$ is a unit. By assumption, if $a \in F(A)^{\Can}_n$ then
\[ ax = a(1 + x) - a \in F(A)^{\Can}_n, \quad xa = (1 + x)a - a \in F(A)^{\Can}_n \]
Since every element of $\JRad(A)$ is nilpotent, this shows that $\JRad(A) \cdot F(A)^{\Can}_n$ and $F(A)^{\Can}_n \cdot \JRad(A)$ are contained in $F(A)^{\Can}_n$. Since $\basefield$ is perfect, $S := A / \JRad(A)$ is separable, so by the Wedderburn-Malt'sev theorem we can write $A \cong S \oplus \JRad(A)$ (as a $\basefield$-vector space). Thus to show the claim it suffices to show that the canonical filtration is closed under multiplying by elements of $S$. After base changing to the algebraic closure, we can assume that $\basefield$ is algebraicially closed and that $S$ is a product of matrix algebras over $\basefield$. The observation above shows that $F(A)^{\Can}_n$ is invariant under multiplication by nilpotent elements of $S$, so it suffices to show the filtration is invariant under multiplication by idempotents of $S$. Let $1 = e_1 + \ldots + e_r$ be a decomposition of $1$ into idempotents, and choose $t \neq 1 \in \basefield$. For any $a \in F(A)^{\Can}_n$, since $e_1 + \ldots + t e_i + \ldots + e_r$ is a unit we have
\[ (1 - t) e_i a = (e_1 + \ldots + e_i + \ldots + e_r) a - (e_1 + \ldots + t e_i + \ldots + e_r) a \in F(A)^{\Can}_n \]
Since $1 - t \neq 0$, this implies that $e_i a \in F(A)^{\Can}_n$. A similar computation shows that $a e_i \in F(A)^{\Can}_n$ as well, so each $F(A)^{\Can}_n$ is invariant under multiplication by idempotents of $S$, which concludes the proof.
\end{proof}

We will now give a result which allows us to detect the canonical filtration in the special case that the associated graded pieces of a filtration happen to be one-dimensional. First, we recall that the rational flag complex $\Delta(\GL(A))$ is the set of one-parameter subgroups of $\GL(A)$ subject to the equivalence relation $\lambda_1(t) \equiv \lambda_2(t)$ if
\[ \lambda_1(t^n) = g \cdot \lambda_2(t^m) \cdot g^{-1} \]
where $n,m \in \bZ$ and $g \in P(\lambda_2)$, with $P(\lambda_2)$ the parabolic group associated to $\lambda_2(t)$. We denote the equivalence class of a one-parameter subgroup $\lambda$ by $\Delta(\lambda)$. See \cite[Chapter 2, \S 2]{mumford1994geometric} for details. Importantly, all elements of an equivalence class define the same filtration.

\begin{lemma} \label{lemma: rational flag complex equivalence preserves algebra filtrations}
Let $\delta(t)$ be a one-parameter subgroup defining an algebra filtration of $A$. Any $\lambda(t) \in \Delta(\delta)$ defines an algebra filtration of $A$, and all these filtrations are the same up to rescaling the weights. In particular, the numerical invariant $\numinvariant(A,-)$ takes the same value on all elements of $\Delta(\delta)$.
\end{lemma}

\begin{proof}
First we show that the every element of the class $\Delta(\delta(t))$ defines the same algebra filtration up to scaling. If $\lambda(t) = \delta(t^m)$ then for all $i$ we have $F(\lambda)_{i} = F(\delta)_{mi}$, so if $\delta(t)$ defines an algebra filtration, then does $\lambda(t)$. From this it is immediate that if $\lambda(t^n) = \delta(t^m)$ and $\delta(t)$ defines an algebra filtration, then so will $\lambda(t)$. It is therefore enough to show that if $\lambda(t) = g \cdot \delta(t) \cdot g^{-1}$ with $g \in P(\delta)$ and $\delta(t)$ defines an algebra filtration, then so does $\lambda(t)$. Observe that if $A_m$ is the weight space for $\delta(t)$ of weight $m$, then $g \cdot A_m$ is the weight space of weight $m$ for $\lambda(t)$. It follows from the construction of the filtration associated to a one-parameter subgroup that $F(\lambda)_n = g \cdot F(\delta)_n$ for all $n$. The parabolic subgroup $P(\delta)$ is exactly the elements of $\GL(A)$ which preserve the filtration $F(\delta)_\bullet$, so we have $g \cdot F(\delta)_n = F(\delta)_n$. It is now clear that $\lambda(t)$ defines an algebra filtration of $A$. Combining all the cases above, we find see that every element of $\Delta(\delta(t))$ defines the same algebra filtration up to rescaling the weights. Since the numerical invariant is scale invariant, this gives the last claim.
\end{proof}

Proposition (\ref{prop: algebra filtration from weighted filtration}) shows that if we fix an (unweighted) vector space filtration $F_{(\bullet)}$ of an algebra $A$, which we write as
\[ A = F_{(0)} \supseteq F_{(1)} \supseteq \ldots \supseteq F_{(p)} \supseteq F_{(p + 1)} = 0 \]
then there is a cone 
\[ \cC(F_{(\bullet)}) \subseteq \{ (w_1, \ldots, w_p) \in \bZ^p \mid w_1 \leq \ldots \leq w_p \} \]
of weight vectors $\bm{w} \in \cC(F_{(\bullet)})$ which define an algebra filtration. Note that if $\bm{w}$ lies in the interior of $F_{(\bullet)}$, so that $w_i < w_{i+1}$ for all $i$, and we assume that $\dim(F_{(i)} / F_{(i + 1)}) = 1$, then the associated graded algebra of the corresponding algebra filtration will have one dimensional graded pieces. On this cone, we can think of the numerical invariant $\numinvariant(A, (F_{(\bullet)}, \bm{w}))$ as a function of the weight vector $\bm{w}$. Furthermore, this makes extends to a continuous function for $\bm{w} \in \cC(F_{(\bullet)})_{\bR} := \cC(F_{(\bullet)}) \otimes \bR^d$. After identifying the set of one-parameter subgroups $\Gamma(S)$ of a rank $n$ torus $S$ with $\bZ^n$, we can also identify $\cC(F_{(\bullet)})$ with a subset of $\Gamma(S)$.

\begin{prop}\label{prop: filtrations with one dimensional graded pieces}
Let $F_{(\bullet)}$ be an unweighted filtration of an unstable algebra $A$
\[ A = F_{(0)} \supseteq F_{(1)} \supseteq \ldots \supseteq F_{(n)} \supseteq F_{(n + 1)} = 0 \]
Let $\bm{w}^\ast$ be the weight vector in $\cC(F_{(\bullet)})_{\bR}$ which minimizes $\numinvariant(A, (F_{(\bullet)}, \bm{w}))$ on $\cC(F_{(\bullet)})_{\bR}$. Suppose that $(F_{(\bullet)}, \bm{w}^\ast)$ is a destabilizing filtration, $w^\ast_i < w^\ast_{i+1}$ for all $i$, and the associated graded algebra has one-dimensional graded pieces. Then $(F_{(\bullet)}, \bm{w}^\ast)$ is the canonical filtration.
\end{prop}

\begin{proof}
Define the following three subsets of $\Delta(\GL(A))$:
\[ \cL(A) := \text{The set of one-parameter subgroups which define algebra filtrations of } A \]
\[ \cN(A) := \text{The set where the function $\nu(A, -)$ is negative in } \Delta(\GL(A)) \]
\[ \cW(S) := \text{The image of } \{ (w_1, \ldots, w_p) \in \bZ^p \mid w_1 \leq \ldots \leq w_p \} \subseteq \Gamma(S) \text{ in } \Delta(S) \]
Thus, $\cN(A) \cap \cW(S) \cap \cL(A)$ is the set of all classes of destabilizing one-parameter subgroups of $\Delta(S)$ which correspond to algebra filtrations.

Let $\lambda_{\Can}$ denote a one-parameter subgroup defining the canonical filtration, and $\lambda^\ast$ denote a one-parameter subgroup corresponding to the weighted filtration in the statement of the proposition (using propositions (\ref{prop: one-parameter subgroups are filtrations of algebras}) and (\ref{prop: algebra filtration from weighted filtration})). Since $\cN(A)_{\bR}$ is convex, there is a line segment $\ell : [0,1] \to \cN(A)_{\bR} \subseteq \Delta(\GL(A))_{\bR}$ joining $\Delta(\lambda^\ast)$ and $\Delta(\lambda_{\Can})$. By \cite[Chapter 2, Lemma 2.9]{mumford1994geometric}, we can choose commuting representatives $\delta^\ast$ and $\delta_{\Can}$ of the classes $\Delta(\lambda^\ast)$ and $\Delta(\lambda_{\Can})$, respectively, and write (in additive notation)
\[ \ell(s) = s \Delta(\delta^\ast) + (1 - s) \Delta(\delta_{\Can}) \]
Since the one-parameter subgroups $\delta^\ast$ and $\delta_{\Can}$ commute and define algebra filtrations by virtue of lemma (\ref{lemma: rational flag complex equivalence preserves algebra filtrations}), it follows that each $\ell(s)$ defines an algebra filtration of $A$. Therefore, $\ell(s) \in [\cN(A) \cap \cL(A)]_{\bR}$ for all $s \in [0,1]$.

Let $T$ denote a maximal torus of $\GL(A)$ determined by a basis compatible with the filtration from the proposition statement. The assumption on the minimizing weights $w^\ast_i$ implies that $\Delta(\lambda^\ast)$ lies in the interior of $\cW(T)$. Thus there is some $s_0 \in (0,1)$ such that $\ell(s_0) \in [\cN(A) \cap \cW(T) \cap \cL(A)]_{\bR}$. Suppose that $\Delta(\lambda^\ast) \neq \Delta(\lambda_{\Can})$. The function $\numinvariant(A, -)$ is strictly convex on $\cN(A)_{\bR}$ (see \cite[Pg. 64]{mumford1994geometric}) which implies that
\begin{align*}
\numinvariant(A, \lambda^\ast) & \leq \numinvariant(A, \ell(s_0)) \\
&< s_0 \numinvariant(A, \lambda^\ast) + (1-s_0) \numinvariant(A, \lambda_{\Can}) \\
& \leq s_0 \numinvariant(A, \lambda^\ast) + (1-s_0) \numinvariant(A, \lambda^\ast) = \numinvariant(A, \lambda^\ast)
\end{align*}
This is a contradiction. Thus we must have $\Delta(\lambda^\ast) = \Delta(\lambda_{\Can})$, which implies that $\lambda^\ast$ corresponds to the canonical filtration.
\end{proof}

\section{Analysis of Stability} \label{section: analysis of stability}

In this section, we give algebraic characterizations of semistable Lie algebras and associative algebras. In the case of Lie algebras, we also show that one can reduce to determining the canonical filtrations of nilpotent Lie algebras.

\subsection{Lie Algebras}

Since $\basefield$ is not necessarily algebraically closed nor characteristic zero, some care must be taken when defining semisimple Lie algebras since the usual equivalent characterizations in the case $\basefield = \bC$ may not hold. We will take the following definition.

\begin{defn} \label{defn: semisimple Lie algebra}
A Lie algebra $L$ is \emph{semisimple} if $\LieRad L = 0$, where $\LieRad L$ is the largest solvable ideal of $L$.
\end{defn}

Note that it is still the case with definition (\ref{defn: semisimple Lie algebra}), $L$ is semisimple if and only if $L$ has no non-zero Abelian ideals; see \cite[\S 5.1]{HumphreysLieAlgebras}.

\begin{thm} \label{thm: A Lie algebra is semistable iff semisimple}
A Lie algebra $L$ is semistable if and only if it is semisimple.
\end{thm}

\begin{proof}
If $L$ is not semisimple then it has an Abelian ideal $I \neq 0$. For instance, one can take $I$ to be the last non-zero term in the derived series of $\LieRad(L)$. The filtration $L \supseteq I$ where $L$ has weight $0$ and $I$ has weight $1$ is then destabilizing. Now assume that $L$ is semisimple. Since each $F(L)^{\Can}_n$ is a Lie ideal by corollary (\ref{coroll: The canonical filtration of a Lie algebra is an ideal filtration}), $I := F(L)_1^{\Can}$ is an ideal of $L$. Let $I^{(n)}$ denote the $n$-th term in the derived series of $I = I^{(0)}$. The compatibility of $F(L)^{\Can}_\bullet$ with the bracket of $L$ implies that $I^{(n)} \subseteq F(L)^{\Can}_{2n}$, hence $I^{(n)} = 0$ for $n \gg 0$ since $F(L)^{\Can}_{2n} = 0$ for $n \gg 0$. This shows that $I$ is solvable, hence $I \subseteq \LieRad L = 0$ and therefore $F(L)^{\Can}_{n} = 0$ for $n > 0$. Thus $- \wt(L, F(L)^{\Can}_\bullet) \geq 0$ and $L$ is semistable.
\end{proof}

\begin{remark}
This generalizes \cite[Theorem 4.3]{lauret2003moment} to perfect fields.
\end{remark}

The following theorem allows us to reduce the computation of the canonical filtration of a general Lie algebra to computing the canonical filtrations of nilpotent or solvable Lie algebras.

\begin{thm} \label{thm: sufficient to compute canonical filtration in the nilpotent case}
Let $L$ be a Lie algebra, and let $\LieRad(L)$ and $\LieNilRad(L)$ denote the radical and nilradical of $L$, respectively. Then 
\[
F(L)^{\Can}_n = 
\begin{cases}
    L & n \leq 0 \\
    F(\LieNilRad(L))^{\Can}_n = F(\LieRad(L))^{\Can}_n & n > 0
\end{cases}
\]
\end{thm}

\begin{proof}
Let $F_\bullet$ be the filtration defined as 
\[
F_n = 
\begin{cases}
    L & n \leq 0 \\
    F(\LieNilRad(L))^{\Can}_n & n > 0
\end{cases}
\]
We claim that $F_\bullet$ is indeed a filtration of Lie algebras. For any $x \in L$, $\ad_{x}$ is a derivation of $\LieNilRad(L)$ since the nilradical is an ideal. Thus $[F_j, F_i] \subseteq F_i$ for all $j$ when $i \leq 0$ by proposition (\ref{prop: action of automorphisms and derivations on canonical filtrations}). When $0 < i \leq j$ $[F_i, F_j] \subseteq F_{i + j}$ is satisfied since $F(\LieNilRad(L))^{\Can}_\bullet$ is a filtration of Lie algebras. Since the Lie bracket is anticommutative, we can assume without loss of generality that $i \leq j$, so the above cases exhaust all possibilities. This proves the claim.

Note that for $n > 0$, $F(L)^{\Can}_n$ is a nilpotent ideal of $L$ and hence is contained in $\LieNilRad(L)$. Define $G_\bullet$ by
\[
G_n = 
\begin{cases}
    \LieNilRad(L) & n \leq 0 \\
    F(L)^{\Can}_n & n > 0
\end{cases}
\]
An argument similar to the one in the previous paragraph shows that $G_\bullet$ is indeed a Lie algebra filtration of $\LieNilRad(L)$. It follows from the definitions that we have
\[ \wt(L, F_\bullet) = \wt(\LieNilRad(L), F(\LieNilRad(L))^{\Can}_\bullet), \quad || F_\bullet || = || F(\LieNilRad(L))^{\Can}_\bullet|| \]
\[ \wt(\LieNilRad(L), G_\bullet) = \wt(L, F(L)^{\Can}_\bullet), \quad || G_\bullet || = || F(L)^{\Can}_\bullet|| \]
These equalities imply that
\[
\nu(L, F(L)^{\Can}_\bullet) = \nu(\LieNilRad(L), G_\bullet) \geq \nu(\LieNilRad(L), F(\LieNilRad(L))^{\Can}_\bullet) = \nu(L, F_\bullet)
\]
This shows that $F(L)^{\Can}_\bullet = F_\bullet$.

To show that $F(L)^{\Can}_n = F(\LieRad(L))^{\Can}$ for $n > 0$, we can repeat the argument above with $\LieRad(L)$ instead of $\LieNilRad(L)$. Indeed, since $\LieRad(L)$ is an ideal, we have $\ad_x \in \Der(\LieRad(L))$. We also observed that $F(L)^{\Can}_n$ is a nilpotent ideal of $L$ for $n > 0$, hence it is solvable and contained in $\LieRad(L)$. The rest of the argument goes through from these observations.
\end{proof}

We say that a Lie algebra $L$ is reductive if $\LieRad(L)$ is the center of $L$. Theorem (\ref{thm: sufficient to compute canonical filtration in the nilpotent case}) allows us to explicitly compute the canonical filtration of reductive Lie algebras.

\begin{coroll} \label{coroll: canonical filtration of a reductive Lie algebra}.
Let $L$ be a reductive Lie algebra. Then the canonical filtration is
\[ F(L)^{\Can}_n = \begin{cases}
  L  & n \leq 0 \\
  \LieRad(L) & n = 1 \\
  0 & n > 1
\end{cases}
\]
\end{coroll}

\begin{proof}
By the theorem (\ref{thm: sufficient to compute canonical filtration in the nilpotent case}), it is enough to determine the canonical filtration of $\LieRad(L)$. Define $F_\bullet$ to be
\[
F_n = 
\begin{cases}
\LieRad(L) & n \leq 1 \\
0 & n > 1
\end{cases}
\]
Observe that $F_\bullet$ is a filtration of Lie algebras since $\LieRad(L)$ is an Abelian ideal, and that
\[ \wt(\LieRad(L), F_\bullet) = ||F_\bullet || = \dim(\LieRad(L)) \]
The Cauchy-Schwarz inequality implies that
\begin{align*}
\numinvariant(\LieRad(L), F(\LieRad(L))^{\Can}_\bullet) &= \frac{- \wt(\LieRad(L), F(\LieRad(L))^{\Can}_\bullet)}{|| F(\LieRad(L))^{\Can}_\bullet||} \\
&\geq - \sqrt{\dim(\LieRad(L))} \\
&= \numinvariant(\LieRad(L), F_\bullet)
\end{align*}
Therefore $F(\LieRad(L))^{\Can}_\bullet = F_\bullet$, which shows that $F(L)^{\Can}_\bullet$ is the claimed filtration.
\end{proof}

\subsection{Associative Algebras}

Let us recall some of the theory of Jacobson radicals of not necessarily unital rings; see \cite[\S 4, Exercies 1-7]{lam1991first}. Let $R$ be a ring. An element $x \in R$ is said to be left (resp. right) quasi-regular if there is an element $y \in R$ such that $y + x - yx = 0$ (resp. $x + y - xy = 0$). The Jacobson radical of $R$ is then defined as
\[ J(R) = \{ x \in R \mid Rx \text{ is left quasi-regular}\} \]
Every element of $\JRad(R)$ is both left and right quasi-regular. If $R$ is unital then $x \in R$ is left (resp. right) quasi-regular if and only if $1 - x$ has a left (resp. right) inverse. In the unital case, $\JRad(R)$ can be characterized as the ideal such that $1-rx$ has a left inverse, which occurs if and only if $rx$ is left quasi-regular. One can also show that, as in the unital case, $\JRad(R)$ is the intersection of the annihilators of all simple left $R$-modules \cite[Exercise 4.6]{lam1991first}. From this description, we see that if $A$ is a $\basefield$-algebra, then in fact $\JRad(A)$ is in fact an \emph{algebra} ideal i.e. is a subspace of $A$: if $c \in \basefield$, $x \in \JRad(A)$, and $m \in M$ for a simple left module $M$ then $(cx)m = c(xm) = 0$, hence $cx \in \JRad(A)$ as well.

\begin{defn} \label{defn: semisimple associative algebra}
A (not necessarily unital) associative algebra $A$ is semisimple if $\JRad(A) = 0$.
\end{defn}

For a (possibly non-unital) associative algebra $A$, the unitization of $A$ is defined to be the algebra $A^\# := \basefield \rtimes A$. As a vector space, $A^{\#}$ is the direct sum $\basefield \oplus A$, and the operations are component-wise addition and multiplication given by
\[ (c_1, a_1) \cdot (c_2, a_2) = (c_1 c_2, c_1 a_2 + c_2 a_1 + a_1 a_2) \]
This makes $A^\#$ into a unital algebra with unit $(1,0)$.

\begin{lemma} \label{lemma: properties of the jacobson radical and semisimple algebras}
Any semisimple associative algebra is unital.
\end{lemma}

\begin{proof}
Let $A$ be a semisimple algebra, possibly non-unital. Consider $(c, x) \in \JRad(A^{\#})$. Since $(c, x)$ is nilpotent, we must have $c = 0$. We will show that $x \in \JRad(A)$, whence we conclude that $x = 0$. For any $y \in A$, we know that $(1,0) - (0,y)(0,x) = (1,-yx)$ must have a left inverse $(t, z)$. It is immediate that we must have $t = 1$. Furthermore,
\[ (1, z) (1,-yx) = (1, z + (-yx) + z(-yx)) = (1,0) \]
which implies that 
\[ (-z) + yx - (-z)yx = 0 \]
which shows that $yx$ is left quasi-regular for any $y \in A$. Therefore $x \in \JRad(A)$, hence $x = 0$ since $A$ is semisimple. This shows that $A^{\#}$ is semisimple and hence is isomorphic to the product of matrix algebras over division rings by the Artin-Wedderburn Theorem. Now $A$ is an ideal of $A^{\#}$, but the only ideals of a product of matrix algebras over division rings are subproducts of these algebras. Thus $A$ is also a product of matrix algebras, and in particular is unital.
\end{proof}

\begin{thm} \label{thm: an associative algebra is semistable iff semisimple}
An associative algebra is semistable if and only if it is semisimple.
\end{thm}

\begin{proof}
Suppose that $A$ is semisimple. By lemma (\ref{lemma: properties of the jacobson radical and semisimple algebras}) it is unital and the Artin-Wedderburn theorem tells us that $A$ is a direct sum of matrix algebras $A_i$, each over some division ring $D_i$ over $\basefield$. Proposition (\ref{prop: canonical filtration for direct sum of algebras}) shows that the canonical filtration of $A$ is a direct sum of the canonical filtrations of the $A_i$, with the weights of the canonical filtration of $A_i$ appropriately scaled. It is therefore enough to show that every matrix algebra is semistable, so we reduce to the case when $A$ is a matrix algebra over a division ring $D$ over $\basefield$.

First we assume that $\basefield$ is algebraically closed, so that $A \cong \Mat_{m}(\basefield)$ for some $m$. Note that if $n > 0$, every element of $F(A)^{\Can}_n$ must be nilpotent. Assume that there is some $x \neq 0 \in F(A)^{\Can}_1$. We claim that this implies that there is some element of $F(A)^{\Can}_1$ which is not nilpotent, which will be a contradiction. By proposition (\ref{prop: action of automorphisms and derivations on canonical filtrations}), $F(A)^{\Can}_1$ is fixed by inner automorphisms i.e. conjugation by units of $A$. This shows that the Jordan canonical form of $x$ is in $F(A)^{\Can}_1$, so we can replace $x$ with its Jordan form. Since $x$ is nilpotent, we have
\[ x = \sum_{j=1}^{m-1} c_j E_{j, j + 1} \]
where the $E_{j,j+1}$ are the standard matrix units, $c_j \in \{ 0, 1\}$, and not all $c_j = 0$ (since $x \neq 0$). Choose $i$ so that $c_i = 1$. Since $F(A)^{\Can}_1$ is also stable under derivations, a straightforward calculation shows that $[E_{i + 1,i}, x] = E_{i+1, i+1} - E_{i,i} \in F(A)^{\Can}_1$. It is clear that $E_{i+1, i+1} - E_{i,i}$ is not nilpotent, which shows the claim. Thus $F(A)^{\Can}_1 = 0$, which implies that $F(A)^{\Can}_n = 0$ for $n > 0$ as well. Therefore $- \wt(A, F(A)^{\Can}_\bullet) \geq 0$, so $A$ is semistable.

Conversely, if $A$ is not semisimple then $\JRad(A)$ is a non-zero nilpotent ideal. The weighted $\JRad(A)$-adic filtration
\[ A = \JRad(A)^0 \supseteq \JRad(A)^1 \supseteq \ldots \supseteq \JRad(A)^r \supseteq 0 \]
where $\JRad(A)^i$ has weight $i$ is destabilizing, hence $A$ is unstable.

Now we prove the theorem assuming that $\basefield$ is merely perfect. By lemma (\ref{lemma: rationality of the canonical filtration}), we see that
\[ \dim_{\overline{\basefield}} F(A_{\overline{k}})^{\Can}_n / F(A_{\overline{k}})^{\Can}_{n + 1} = \dim_{\basefield} F(A)^{\Can}_n / F(A)^{\Can}_{n + 1} \]
This implies that
\begin{equation} \label{eqn: invariance of wt under basechange}
\wt(A, F(A)^{\Can}_\bullet) = \wt(A_{\overline{\basefield}}, F(A_{\overline{\basefield}})^{\Can}_\bullet)
\end{equation}
hence $A$ is semistable if and only if $A_{\overline{\basefield}}$ is semistable. Assume $A$ is semisimple. Since $\basefield$ is perfect, $A_{\overline{k}}$ is semisimple \cite[Corollary 6.1.4]{drozd2012finite} and the argument above in the algebraically closed case shows that $A_{\overline{k}}$ is semistable. Hence $A$ is semistable too. If $A$ is semistable, we claim that $A_{\overline{\basefield}}$ is also semistable. Indeed, since $\basefield$ is perfect any optimal destabilizing one-parameter subgroup $\overline{\lambda}$ of $A_{\overline{k}}$ is defined over $\basefield$, so that there exists a one-parameter subgroup $\lambda$ defined over $\basefield$ such that $\lambda \times_{\Spec(k)} \Spec(\overline{k}) = \overline{\lambda}$. From the equality (\ref{eqn: invariance of wt under basechange}), we see that $\lambda$ would destabilize $A$, proving the claim. The argument above then shows that $A_{\overline{\basefield}}$ is semisimple. Since algebraically closed fields are perfect, $A_{\overline{\basefield}}$ is separable. By \cite[Corollary 6.1.7]{drozd2012finite}, $A$ is separable and in particular semisimple.
\end{proof}

\begin{remark}
This generalizes \cite[Theorem 4.6]{zhang2023moment} to perfect fields.
\end{remark}

As an application of the above theorem, we show how to compute the canonical filtration of some algebras constructed from semisimple algebras. For an associative algebra $A$ and a $A$-bimodule $M$, recall that the trivial extension of $A$ by $M$ is the algebra which is equal to $A \oplus M$ as a vector space with multiplication given by
\[ (a, m) (a', m') = (a a', a m' + m a') \]
We denote this algebra by $A \rtimes_0 M$, where the $0$ subscript indicates that we are extending the multiplication of $A$ by the $0$-cocycle.

\begin{prop}[Trivial Extensions of Semisimple Algebras]
Let $A$ be a semisimple algebra and $M$ an $A$-bimodule. Then the canonical filtration of the trivial extension $A \rtimes M$ is
\[ F(A \rtimes_0 M)^{\Can}_n = \begin{cases}
  A \rtimes_0 M  & n \leq 0 \\
  M & n = 1 \\
  0 & n > 1
\end{cases}
\]
\end{prop}

\begin{proof}
Let $H_\bullet$ denote the filtration defined in the statement of the proposition. One can readily check that $\delta_M : A \rtimes_0 M \to M$, $\delta_M((a,m)) = (0, m)$ is a derivation. Proposition (\ref{prop: action of automorphisms and derivations on canonical filtrations}) then implies that as vector spaces
\[ F(A \rtimes_0 M)^{\Can}_n = (F(A \rtimes_0 M)^{\Can}_n\cap A) \oplus (F(A \rtimes_0 M)^{\Can}_n \cap M) \]
Define the filtrations $F_\bullet = F(A \rtimes_0 M)^{\Can}_\bullet \cap A$ and $G_\bullet = F(A \rtimes_0 M)^{\Can}_\bullet \cap M$ of $A$ and $M$ as vector spaces, respectively. Set
\[ L_A =\sum_{n \in \bZ} n \dim(F_ n/ F_{n+1}), \quad B_A =\sum_{n \in \bZ} n^2 \dim(F_ n/ F_{n+1}) \]
\[ L_M =\sum_{n \in \bZ} n \dim(G_ n/ G_{n+1}), \quad B_M =\sum_{n \in \bZ} n^2 \dim(G_ n/ G_{n+1}) \]
and note that because of the direct sum decomposition above,
\[ \nu(A \rtimes M, F(A \rtimes_0 M)^{\Can}_\bullet) = - \frac{L_A + L_M}{\sqrt{B_A +B_M}}  \]
Since $A$ is semisimple, $L_A \leq 0$ by theorem (\ref{thm: an associative algebra is semistable iff semisimple}), and we always have $B_A \geq 0$. Therefore,
\[ \nu(A \rtimes_0 M, F^{\Can}_{\bullet}) \geq - \frac{L_M}{\sqrt{B_M}} \geq - \frac{\dim(M)}{\sqrt{\dim(M)}} = \nu(A \rtimes_0 M, H_\bullet) \]
where the right-most inequality is due to the Cauchy-Schwarz inequality. This implies that $H_\bullet = F(A \rtimes_0 M)^{\Can}_\bullet$.
\end{proof}

Let $A$ and $B$ be associative algebras and $M$ and $(A,B)$-bimodule. The triangular algebra $T(A, M, B)$ is defined as
\[ T(A, M, B) = \left\{
\begin{bmatrix}
a & m \\
0 & b \\
\end{bmatrix} \colon a \in A, b \in B, m \in M \right\}
\]
with multiplication given by
\[ \begin{bmatrix}
a & m \\
0 & b \\
\end{bmatrix} \begin{bmatrix}
a' & m' \\
0 & b' \\
\end{bmatrix}
= 
\begin{bmatrix}
a a' & a m' + m b' \\
0 & b b' \\
\end{bmatrix}
\]

\begin{prop}
Let $A$ and $B$ be semisimple algebras and $M \neq 0$ an $(A, B)$-bimodule. Then the canonical filtration of $T := T(A, M, B)$ is given by
\[ F(T)^{\Can}_n = \begin{cases}
  T  & n \leq 0 \\
  M & n = 1 \\
  0 & n > 1
\end{cases}
\]
\end{prop}

\begin{proof}
First observe that 
\[ I = \left\{
\begin{bmatrix}
0 & m \\
0 & 0 \\
\end{bmatrix} \colon m \in M \right\} \]
is an ideal of $T$, and that $T / I \cong A \oplus B$ is semisimple. Set
\[ F_\bullet = F(T)^{\Can}_\bullet \cap I, \quad G_\bullet = F(T)^{\Can}_\bullet / (F(T)^{\Can}_\bullet \cap I) \]
Via the isomorphisms $F(T)^{\Can}_n / (F(T)^{\Can}_n \cap I) \cong (F(T)^{\Can}_n + I) / I$, a routine check shows that we can view $G_\bullet$ as an algebra filtration of $A \oplus B$. Set
\[ L_{T / I} = \sum_{n \in \bZ} n \dim(G_n / G_{n+1}), \quad B_{T/ I} = \sum_{n \in \bZ} n^2 \dim(G_n / G_{n+1}) \]
\[ L_{I} = \sum_{n \in \bZ} n \dim(F_n / F_{n+1}), \quad B_{I} = \sum_{n \in \bZ} n^2 \dim(F_n / F_{n+1}) \]
We have
\[ \wt(T, F(T)^{\Can}_\bullet) = L_{T / I} + L_I, \quad || F(T)^{\Can}_\bullet ||^2 = B_{T / I} + B_I \]
Let $H_\bullet$ denote the filtration from the proposition statement. Since $L_{T / I} \leq 0$ and $B_{T / I} \geq 0$, we see that
\[ \nu(T, F(T)^{\Can}_\bullet) = - \frac{L_I + L_{T/I}}{\sqrt{B_I + B_{T / I}}} \geq - \frac{L_I}{\sqrt{B_I}} \geq - \sqrt{\dim(M)} = \nu(T, H_\bullet) \]
Where the first inequality follows from the semistability of $T/I$ and the second inequality follows from the Cauchy-Schawrz inequality. This shows that $H_\bullet = F(T)^{\Can}_\bullet$ as claimed.
\end{proof}

\section{A Strategy to Compute Canonical Filtrations} \label{section: strategy to compute canonical filtrations}

In general, computing the canonical filtration of an algebra seems to be a difficult problem. Even if we know the unweighted filtration associated to the canonical filtration, determining the optimal weights is a non-trivial task. In this section, we outline a strategy to reduce the computation of the canonical filtration of an algebra $A$ to a smaller set of optimization problems when $A$ has a sufficiently large set of automorphisms and/or derivations. In some instances, such as when $A$ is $\bZ^r$-graded with one dimensional homogeneous components, this set will even be finite. A note on conventions: bold face text (e.g. $\bm{x}$) denotes vectors, $\bm{0}$ denotes the zero vector, $C^t$ denotes the transpose of a matrix $C$, and $\bm{x} \geq \bm{y}$ means that each coordinate of $\bm{x}$ is greater than the corresponding coordinate of $\bm{y}$.

\subsection{Some Convex Optimization}
Fix a vector $\bm{d} \in \bZ_{> 0}^{r + 1}$ and define an inner product $\langle -, - \rangle_{\bm{d}}$ on $\bR^{r + 1}$ as
\[ \langle \bm{u}, \bm{v} \rangle_{\bm{d}} = \sum_{i = 0}^r u_i v_i d_i \]
This induces a norm on $\bR^{r + 1}$, which we denote by $|| - ||_{\bm{d}}$. Define a function $f_{\bm{d}}: \bR^{r+1} \setminus \bm{0} \to \bR$ by
\[ f_{\bm{d}}(\bm{x}) =\frac{\langle \bm{1}, \bm{x} \rangle_{\bm{d}}}{|| \bm{x} ||_{\bm{d}}} \]
where $\bm{1} \in \bR^{r + 1}$ is the vector all of whose entries are equal to $1$. Note that our numerical invariant $\nu$ is the negative of such a function. Indeed, suppose we have a filtration $F_\bullet$ corresponding to an unweighted filtration
\[ A = F_{(0)} \supsetneq F_{(1)} \supsetneq \ldots \supsetneq F_{(r)} \supsetneq 0 \]
and a weighting given by a vector $\bm{w} \in \bR^{r+1}$. If we set $d_i = \dim(F_{(i)} / F_{(i + 1)} )$, we see that
\[ \wt(A, F_\bullet) = \langle \bm{1}, \bm{w} \rangle_{\bm{d}}, \quad ||F_\bullet|| = || \bm{w} ||^2_{\bm{d}} \]
and so $\nu(A, F_\bullet) = - f_{\bm{d}}(\bm{w})$. We can also define a filtration of $A$ by choosing a basis $\cB = \{f_1, \ldots, f_n \}$, weights $w_i$ corresponding to $f_i$, and then extending this to a weight function on $A$ corresponding to a filtration $F_\bullet$. In this case we have
\[ \wt(A, F_\bullet) = \sum_{i=1}^n w_i = \langle \bm{1}, \bm{w} \rangle_{\bm{1}}, \quad ||F_\bullet|| = \sum_{i=1}^n w_i^2 = || \bm{w} ||^2_{\bm{1}} \]
where we take $\bm{w}^{\Transpose} = (w_1, \ldots, w_n)$. In this scenario, the numerical invariant is $\numinvariant(A, F(w)_\bullet) = -f_{\bm{1}}(\bm{x})$. In any case, to minimize, $\numinvariant$ we wish to \emph{maximize} $f_{\bm{d}}(\bm{x})$ over some convex set $S$. For the purposes of computing canonical filtrations, it will be convenient to consider the function
\[ g_{\bm{d}}(\bm{x}) = \frac{1}{2} || \bm{x} - \bm{1} ||^2_{\bm{d}} \]
whose optimizer is closely related to the optimizer of $f_{\bm{d}}$.

\begin{prop} \label{prop: equivalence of optimization models}
Let $f_{\bm{d}}(\bm{x})$ and $g_{\bm{d}}$ be as above. Set
\[ S = \{ \bm{x} \in \bR^d \mid \langle u, \bm{1} \rangle_{\bm{d}} \geq 0 \} \]
which is a closed convex set. Then we have
\begin{enumerate}
    \item[(i)] If a point $\bm{x}^\ast \in S$ is a global minimum for $g_{\bm{d}}|_S$ then it is a global maximum for $f_{\bm{d}}|_S$.
    \item[(ii)] If a point $\bm{y}^\ast \in S$ is a global maximum for $f_{\bm{d}}|_S$ then
    \[ \bm{x}^\ast = \frac{\langle \bm{1}, \bm{y}^\ast \rangle_{\bm{d}}}{|| \bm{y}^\ast ||^2_{\bm{d}}} \cdot \bm{y}^\ast \]
    is a global minimum for $g_{\bm{d}}|_S$.
\end{enumerate}
If $\bm{y}^\ast$ is a maximizer for $f_{\bm{d}}|_S$, we define its normalization to be the point $\bm{x}^\ast$ in (ii).
\end{prop}

\begin{proof}
Endow $\bR^{r+1}$ with the inner product $\langle -, - \rangle_{\bm{d}}$. It is clear that $g_{\bm{d}}|_S$ is minimized by the point in $S$ closest to $\bm{1}$ in the metric induced by $\bm{d}$, and that it is a strictly convex function, hence has a unique minimizer. For any $\bm{u} \in S$, let $P(\bm{u})$ be the orthogonal projection of $\bm{1}$ onto the line through $\bm{u}$, and let $z(\bm{u}) = \bm{1} - P(\bm{u})$. For $\bm{u}, \bm{v} \in \bR^{r+1}$, properties of orthogonal projections show that
\[ || P(\bm{u}) ||^2_{\bm{d}} + || z(\bm{u}) ||^2_{\bm{d}} = || \bm{1} ||^2_{\bm{d}} = || P(\bm{v}) ||^2_{\bm{d}} + || z(\bm{v}) ||^2_{\bm{d}} \]
\begin{equation} \label{eqn: equivalence of optimization models eqn 1}
\implies || P(\bm{u}) ||^2_{\bm{d}} - || P(\bm{v}) ||^2_{\bm{d}} = || z(\bm{v}) ||^2_{\bm{d}} - || z(\bm{u}) ||^2_{\bm{d}} \end{equation}
For $\bm{u}, \bm{v} \in S$, we have
\[ || P(\bm{u}) ||_{\bm{d}} = \left| \left| \frac{\langle \bm{u}, \bm{1} \rangle_{\bm{d}}}{|| \bm{u} ||_{\bm{d}}^2} \bm{u} \right| \right|_{\bm{d}} = \frac{\langle \bm{u}, \bm{1} \rangle_{\bm{d}}}{|| \bm{u} ||_{\bm{d}}^2} || \bm{u} ||_{\bm{d}} = \frac{\langle \bm{u}, \bm{1} \rangle_{\bm{d}}}{|| \bm{u} ||_{\bm{d}}} = f_{\bm{d}}(\bm{u}) \]
where we have used the fact that $| \langle \bm{u}, \bm{1} \rangle_{\bm{d}} | = \langle \bm{u}, \bm{1} \rangle_{\bm{d}}$ since $\bm{u} \in S$. Observe that $|| z(\bm{u}) ||_{\bm{d}}$ is the distance between $\bm{1}$ and the line spanned by $\bm{u}$, and that $||z(\bm{u})||_{\bm{d}}$ is minimized on $S$ precisely when $g_{\bm{d}}|_S$ is minimized. Thus, equation (\ref{eqn: equivalence of optimization models eqn 1}) shows that $f_{\bm{d}}|_S$ is maximized by the unique minimizer of $g_{\bm{d}}|_S$, which gives (i). Conversely, if $\bm{y}^\ast$ is a maximizer of $f_{\bm{d}}|_S$ then so is any point on the line spanned by $\bm{y}^\ast$ since $f_{\bm{d}}|_S$ is scale invariant. In particular, the projection $P(\bm{y}^\ast)$ is also a maximizer for $f_{\bm{d}}|_S$, which also minimizes $g_{\bm{d}}|_S$ thanks to equation (\ref{eqn: equivalence of optimization models eqn 1}). The definition of $\bm{x}^\ast$ in (ii) is precisely $P(\bm{y}^\ast)$, proving (ii).
\end{proof}

Since $g_{\bm{d}}$ is a convex function, to check that $\bm{x}^\ast$ is a minimizer it is sufficient to check that $\bm{x}^\ast$ satisfies the Karush-Khun-Tucker (KKT) conditions. We state a slightly simplified version of the KKT conditions (without the equality constraints) for reference.

\begin{thm}[The KKT Conditions] \label{thm: KKT conditions}
Let $C \in \Mat_{m \times n}(\bR)$ and $f : \bR^n \to \bR$ be a convex function. Then $\bm{x}^\ast$ is a minimum for $f$ subject to the constraint $C \bm{x}^\ast \leq \bm{0}$ if and only if the following hold:
\begin{enumerate}
    \item[(i)] $\nabla f(\bm{x}^\ast) + C^{\Transpose} \bm{\lambda} = \bm{0}$, with $\bm{\lambda} \geq \bm{0} \in \bR^m$.
    \item[(ii)] $C \bm{x}^\ast \leq \bm{0}$.
    \item[(iii)] $\bm{\lambda^{\Transpose}} (C \bm{x}^\ast) = \bm{0}$.
\end{enumerate}
\end{thm}

\subsection{The Strategy} \label{subsection: Strategy}

Suppose $\cB = \{ a_1, \ldots, a_d \}$ and $w_1, \ldots, w_d$ is a weighted algebra basis of $A$ as in definition (\ref{defn: weighted algebra basis}). Recall that the weights $w_i$ must satisfy the weight inequalities coming from the multiplication of $A$. That is, if 
\[ \mu(a_i, a_j) = \sum_{m=1}^d c_{ij}^m a_m \]
then we have $w_i + w_j - w_m \leq 0$ whenever $c_{ij}^m \neq 0$. Let $\bm{w} = (w_1, \ldots, w_d)^{\Transpose} \in \bR^d$ and let $C$ be a matrix with one row $(\bm{e}_i + \bm{e}_j - \bm{e}_m)^{\Transpose}$ for each weight inequality as above, where the $\bm{e}_r$ are the standard basis vectors. Then we have $C \bm{w} \leq \bm{0}$. Conversely, given the basis $\cB$ and a vector $\bm{w}$ satisfying $C \bm{w} \leq \bm{0}$, we see that $w_1, \ldots, w_d$ satisfy the weight inequalities and hence define a weighted algebra basis. Thus for a fixed basis $\cB$, we can find the maximum of $\numinvariant(A, F_\bullet)$ over filtrations $F_\bullet$ compatible with $\cB$ by solving the following optimization problem:
\begin{equation} \label{eqn: the convex optimization problem}
\begin{aligned}
\textit{minimize} \quad & g(\bm{w}) = \frac{1}{2} \sum_{i=1}^d (w_i - 1)^2 \\ 
\textit{s. t.} \quad & C \bm{w} \leq \bm{0} \\
\end{aligned}
\end{equation}
Of course, to compute the canonical filtration we must solve this problem over \emph{every} such basis. Thus if we can determine a set $S$ of bases which are compatible with the canonical filtration, we can reduce the computation of the canonical filtration to solving problem (\ref{eqn: the convex optimization problem}) for each basis in $S$. In particular, if we can find a finite set $S$ of bases, we need only solve finitely many optimization problems.

We can use this observation to compute the canonical filtration of $\bZ^r$-graded algebras. By corollary (\ref{coroll: canonical filtration is a direct sum of graded pieces}), an ordered basis $\cB$ consisting of homogeneous elements of $A$ is compatible with the canonical filtration, so we need only solve the optimization problem for bases consisting of homogeneous elements. In particular, if each homogeneous component has dimension one then we can use a single basis, and thus computing the canonical filtration amounts to solving a single instance of problem (\ref{eqn: the convex optimization problem}). If $A$ is graded but the homogeneous components are not all one-dimensional but we have additional knowledge about the automorphism group, we may still be able to reduce to solving a restricted number of instances of the optimization problem, e.g. if there are automorphisms which respect the given grading. Thus it is natural to investigate and compute canonical filtrations for $\bZ^r$-graded algebras, and we will do precisely this in \autoref{section: examples and computations}.

In general, the optimal weights for a given basis $\cB$ need not satisfy all constraints with equality, even we start with a basis consisting of homogeneous elements for some grading. Thus it is necessary to check that $\bm{\lambda} \geq \bm{0}$ in the KKT conditions (\ref{thm: KKT conditions}). However, if a feasible set of weights for $\cB$ satisfy the constraints with equality, so that we are giving a grading on $A$ where the elements of $\cB$ are homogeneous, then checking the KKT conditions is much easier. In this case, (ii) and (iii) of the KKT conditions are automatic, so we need only check (i). Computing the vector $\bm{\lambda}$ in condition (i) may be difficult if the constraint matrix $C$ is complicated, but fortunately (i) can be verified indirectly thanks to the following lemma.

\begin{lemma}[Farkas' Lemma] \label{lemma: Farkas' lemma}
Let $\bm{w} \in \bR^n$ and $C \in \Mat_{m \times n}(\bR)$. Then exactly one of the following is true:
\begin{itemize}
    \item[(i)] There is some $\bm{\lambda} \in \bR^m$ such that $\bm{w} = C^{\Transpose} \bm{\lambda}$ and $\bm{\lambda} \geq 0$.
    \item[(ii)] There is some $\bm{y} \in \bR^n$ such that $C \bm{y} \leq 0$ and $\bm{y}^{\Transpose} \bm{w} > 0$
\end{itemize}
\end{lemma}

See \cite[Proposition 1.8]{Ziegler} for details. We can use this lemma to verify (i) in the KKT conditions in the case where the $C \bm{w^\ast} = \bm{0}$ for our candidate optimizer $\bm{w}^\ast$: assume we have some $\bm{y} = (y_1, \ldots, y_d)^t$ such that $C \bm{y} \leq \bm{0}$, and then we show that
\[ \bm{y}^t ( - \nabla g_{\bm{1}}(\bm{w}^\ast) ) = \sum_{i = 1}^d (1 - w^\ast_i) y_i \leq 0 \] Then (ii) is false in Farkas' Lemma (\ref{lemma: Farkas' lemma}), so condition (i) of the same lemma must hold, hence condition (i) of the KKT conditions is verified.

We conclude this section with a trick which allows us to handle some cases when $\basefield$ has positive characteristic. Note that if the algebraic relations of the elements in the basis $\cB$ do not depend on $\basefield$, then the constraint matrix in the optimization problem above will be identical to the constraint matrix in the characteristic $0$ case. The resulting optimization problem will thus be equivalent to the characteristic $0$ case, thus we can compute the canonical filtration assuming $\basefield$ has characteristic $0$. For example, this is the case for an algebra of the form $\basefield[x_1, \ldots, x_n] / I$ where $I$ is an ideal consisting of unital monomials.

\begin{prin}[The Characteristic $0$ Principle] \label{prin: char 0 principle}
If the relations of a given basis $\cB$ do not depend on $\basefield$, then we may assume $\basefield$ has characteristic $0$ when computing the optimizing weights for $\cB$.
\end{prin}

\section{Examples and Computations} \label{section: examples and computations}

In this section, we illustrate the strategy outlined in \autoref{subsection: Strategy} to compute the canonical filtration in some examples.

\subsection{Monomial Algebras}

Let $A$ be a commutative finite-dimensional algebra of the form
\[ A \cong \basefield[x_1, \ldots, x_n] / I \]
where $I$ is generated by monomials in the $x_i$ with coefficient $1$. We will call an algebra of this form a \emph{monomial algebra}. For commutative rings, every Artinian ring is a direct sum of local Artinian rings, so by proposition (\ref{prop: canonical filtration for direct sum of algebras}), it suffices to consider the case of local monomial algebras. These algebras are well-suited to our techniques since the standard $\bZ^n$-grading is such that the dimension of the degree $\alpha$ elements is at most one.

\begin{prop} \label{prop: the canonical filtration of a monomial algebra is a filtration by ideals}
If $A$ be a commutative monomial algebra then $F(A)^{\Can}_\bullet$ is a filtration of $A$ by ideals.
\end{prop}

\begin{proof}
Since the relations of the monomial basis do not depend on $\basefield$, we can apply the characteristic $0$ principle (\ref{prin: char 0 principle}) and assume that $\basefield$ has characteristic $0$. First we assume that $A$ is local. Let $x^\alpha \neq 1 \in F(A)^{\Can}_m$ be a monomial. For any $i$ such that $x_i$ divides $x^\alpha$, consider the derivation $\nabla_{x_i} = x_i \partial_{x_i}$ of $\basefield[x_1, \ldots, x_n]$. This preserves any monomial ideal, so it descends to a derivation of $A$. For any $f \in A$, $f \nabla_{x_i}$ is a derivation of $A$ and we have $f \nabla_{x_i}(x^\alpha) = \alpha_i f x^\alpha \in F(A)^{\Can}_m$. We can scale away $\alpha_i$ thanks to the characteristic $0$ assumption, which implies that $f x^\alpha \in F(A)^{\Can}_m$. 

Define $G_\bullet$ by setting $G_n = F(A)^{\Can}_n$ if $n > 0$ and $G_n = A$ if $n \leq 0$. This is an algebra filtration of $A$. Indeed, it suffices to check that if $f \notin G_1$ and $g \in G_n$ for $n \geq 1$ then $fg \in G_n$. Since $A$ is local, we can write $g = c + m$ for some $c \in \basefield$ and $m \in \mathfrak{m}_A$, where $\mathfrak{m}_A$ is the maximal ideal of $A$. Observe that $g$ must be nilpotent since $n \geq 1$, hence $c = 0$. Therefore, $g = \sum_{\alpha} c_\alpha x^\alpha$ over all $\alpha \neq 0$, and so $f g \in G_n$ by the preceding paragraph. One checks that $\nu(A, G_\bullet) \leq \nu(A, F(A)^{\Can}_\bullet)$; this follows from the fact $G_\bullet$ collapses all non-positively indexed filtered pieces of $F(A)^{\Can}_\bullet$ into $G_0$. This shows that $G_\bullet = F(A)^{\Can}_\bullet$.

In the non-local case, we know from proposition (\ref{prop: canonical filtration for direct sum of algebras}) that the canonical filtration of a non-local monomial algebra $A = \oplus_{i=1}^s A_i$ is the sum of the filtrations of the $A_i$ with the indices rescaled by positive integers. Since each factor in $A$ is a local algebra, whose canonical filtration is a filtration by ideals, the same is true for $A$.
\end{proof}

We now compute the canonical filtration of some explicit monomial algebras. We remark that since the canonical filtration of a commutative monomial algebra $A$ is a filtration by ideals by the proposition above, proposition (\ref{prop: filtration by ideals}) implies that the unit of $A$ has weight $0$ in all the following computations.

\begin{example} \label{example: jordan algebras}
The canonical filtration of $A = \basefield[x] / (x^n)$ is given by the 
assignment of weights $w(x^i) = i$ for $0 \leq i \leq n$.
\end{example}

\begin{proof}
We will use the quadratic optimization function $g_{\bm{1}}$ from (\ref{prop: equivalence of optimization models}). Since the proposed weights define a grading, we need only show that the normalized weights $\frac{3i}{2n + 1}$ satisfies the condition of Farkas' lemma (\ref{lemma: Farkas' lemma}). To this end, let $y_i$ be sequence of real numbers satisfying the constraints imposed by the monomials $x^i$ for $i \geq 0$. Since the unit has weight $0$, it is immediate that $y_0 \leq 0$ and we ignore this term in our computation. We must verify that
\[ \sum_{i=1}^n \left( 1 - \frac{3i}{2n+1} \right) y_i = \frac{1}{2n+1} \sum_{i=1}^n (2n + 1 - 3i) y_i \leq 0 \]
 We claim that
\[ \frac{1}{2n+1} \sum_{i=1}^n (2n + 1 - 3i) y_i = \frac{1}{2n+1} \sum_{i=1}^n \sum_{j=1}^{i-1} (y_j + y_{i-j} - y_i) \]
For a fixed index $q$, we consider the sum of the coefficients of $y_q$ on the right-hand side and show that it equals $2n + 1 - 3q$. The number of times $y_q$ appears with a coefficient of $-1$ is clearly $q - 1$. For each $i > q$, $y_q$ will appear in the inner sum of the right-hand side exactly twice, and there are $n - q$ such $i$. Therefore, the number of times $y_q$ appears with a coefficient of $+1$ is $2(n - q)$. The coefficient of $y_q$ in the right-hand sum is then $2(n - q) - (q - 1) = 2n + 1 - 3q$ as claimed. Since $y_j + y_{i - j} - y_i \leq 0$ for each pair of $i,j$, the claim follows.
\end{proof}

As a consequence, we can determine the canonical filtrations of commutative algebras $A = \basefield[x] / (f(x))$ where $f(x)$ is monic $f(x)$. Suppose $f(x)$ factors as
\[ f(x) = \prod_{i=1}^m (x - r_i)^{n_i} \]
From the Chinese remainder theorem and change of coordinates $x \mapsto x + r_i$, we see that
\[ A \cong \bigoplus_{i=1}^m \basefield[x] / (x^n_i) \]
so we can compute the canonical filtration using example (\ref{example: jordan algebras}) and proposition (\ref{prop: canonical filtration for direct sum of algebras}).

One might hope that since monomial algebras are graded, their canonical filtrations arise from a grading. The following example, selected for its relative simplicity, shows that this need not be the case.

\begin{example} \label{example: small monomial algebra whose canonical filtration does not arise from a grading}
The canonical filtration $F(A)^{\Can}_\bullet$ of the algebra $A = \basefield[x, y] / (x^4, x^2 y, x y^2, y^4)$ is given by the assignment of weights
\[ w(1) = 0, \quad w(x^i) = w(y^i) = 3i, \quad w(xy) = 7 \]
for $1 \leq i \leq 3$. In particular, $F(A)^{\Can}_\bullet$ does not arise from a grading of $A$ since $w(x) + w(y) < w(xy)$.
\end{example}

\begin{proof}
The normalized weight function $v$ of the proposed weight function $w$ is given by
\[ v(x^i) = \frac{3i}{7}, \quad v(y^i) = \frac{3i}{7}, \quad v(xy) = 1 \]
We verify the KKT conditions. The constraint matrix for the ordered basis 
\[ \{ x, x^2, x^3, y, y^2, y^3, xy \} \]
is 
\[
C = 
\begin{bmatrix}
2 & -1 & 0 & 0 & 0 & 0 & 0 \\
1 & 1 & -1 & 0 & 0 & 0 & 0 \\
0 & 0 & 0 & 2 & -1 & 0 & 0 \\
0 & 0 & 0 & 1 & 1 & -1 & 0 \\
1 & 0 & 0 & 1 & 0 & 0 & -1 \\
\end{bmatrix}
\]
Let $\bm{x}$ be the vector of normalized weights, ordered according to the basis above:
\[ \bm{x}^{\Transpose} = 
\begin{bmatrix}
\frac{3}{7} & \frac{6}{7} & \frac{9}{7} & \frac{3}{7} & \frac{6}{7} & \frac{9}{7} & 1 \\
\end{bmatrix}
\]
One can check that $C \bm{x} \leq \bm{0}$, and if we let
\[ \bm{\lambda}^{\Transpose} =
\begin{bmatrix}
\frac{1}{7} & \frac{2}{7} & \frac{1}{7} & \frac{2}{7} & 0
\end{bmatrix}
\]
then $\bm{1} - \bm{x} = C^{\Transpose} \bm{\lambda}$. For the complimentary slackness condition, one checks that $\bm{\lambda}^{\Transpose} C \bm{x} = \bm{0}$. This shows that all KKT conditions are satisfied.
\end{proof}

Even though the canonical filtration in example (\ref{example: small monomial algebra whose canonical filtration does not arise from a grading}) does not arise from a grading, the corresponding weight function $w$ is homogeneous in the sense that for any monomial $f$ and $m \in \bZ_{\geq 0}$ such that $f^m \neq 0$,  $w(f^m) = m \cdot w(f)$. We next examine a class of monomial algebras where even this property does not hold.

Let $S = \basefield[x_1, \ldots, x_n]$, $\mathfrak{m} = (x_1, \ldots, x_n)$, and set $S_{n,r} = S / \mathfrak{m}^r$. We assume that the characteristic of $\basefield$ is either $0$ or greater than $r$. If $\partial_{x_i}$ denotes the usual differential operator on $S$, observe that $x_j \partial_{x_i}$ preserves the total degree of any monomial. Therefore, this derivation preserves the ideal $\mathfrak{m}^r$ and descends to a derivation of $S_{n,r}$ for any $r$. Since the canonical filtration is stable under derivations, if $F(S_{n,r})^{\Can}_m$ contains one monomial of degree $d$ then one can use the derivations $x_j \partial_i$ to show that it must contain all monomials of degree $d$. Therefore, the canonical filtration of $S_{n, r}$ is the $\mathfrak{m}$-adic filtration with some assignment of weights. Even though the underlying filtration is known, computing $F(S_{n,r})^{\Can}_\bullet$ is still challenging. However, for small values of $r$ we can compute the filtration exactly.

\begin{example} \label{example: S_{n,4}}
When $2 \leq n \leq 4$, the weights of the canonical filtration of $S_{n, 4}$ are given by $w(\mathfrak{m}^i) = i$. When $n > 5$, the weights of positive powers of $\mathfrak{m}$ are given by 
\[ w(\mathfrak{m}) = 2n + 10, \quad w(\mathfrak{m}^2) = n^2 + 3n + 8, \quad w(\mathfrak{m}^3) = n^2 + 5n + 18 \]
\end{example}

\begin{proof}
We will formulate the optimization problem using the function $g_{\bm{d}}(\bm{x})$ as in proposition (\ref{prop: equivalence of optimization models}), where the $d_i$ are the number of monomials of degree $i$:
\[ d_i = \binom{n - 1 + i}{i} \]
Set $\bm{d} = (d_1, d_2, d_3)^{\Transpose}$. The constraint matrix is
\[
C = 
\begin{bmatrix}
2 & -1 & 0 \\
1 & 1 & -1
\end{bmatrix}
\]
We set $z_i = d_i(1 - x_i)$, so that $\nabla g_{\bm{d}}(\bm{x}) = (z_1, z_2, z_3)^{\Transpose}$. Then the KKT conditions require that the optimal $\bm{x}$ satisfies
\[ z_1 = 2\lambda_1 + \lambda_2, \quad z_2 = - \lambda_1 + \lambda_2, \quad z_3 = - \lambda_2  \]
where $\lambda_i \geq 0$. Solving for the $\lambda_i$ in terms of the $z_i$, we find that
\[ \lambda_1 = \frac{1}{3} (z_1 - z_2), \quad \lambda_2 = \frac{1}{3} (z_1 + 2z_2) = - z_3 \]
and the equation for $\lambda_2$ implies that $z_1 + 2 z_2 + 3 z_3 = 0$. The complementary slackness conditions can be written as
\[ (z_1 - z_2)(2x_1 - x_2) = 0, \quad (z_1 + 2 z_2)(x_1 + x_2 - x_3) = 0 \]
We cannot have $z_1 - z_2 = 0$ and $z_1 + 2 z_2 = 0$ at once. Indeed, these two equalities plus the condition that $z_1 + 2 z_2 + 3 z_3 = 0$ would imply that $x_1 = x_2 = x_3 = 1$, but this is not a feasible point. Next we claim that we must have $z_1 + 2 z_2 \neq 0$. Assume to the contrary that $z_1 + 2 z_2 = 0$. By the preceding argument we must have $z_1 - z_2 \neq 0$, hence $\lambda_1 \neq 0$. The equality $z_1 + 2 z_2 = 0$ implies that $\lambda_2 = 0$, hence $z_3 = 0$ so $x_3 = 1$. The complementary slackness conditions imply that $x_2 - 2 x_1 = 0$, but $2 z_2 = - z_1$ implies that
\[ x_2 = \frac{d_1}{2 d_2}(1 - x_1) + 1 = 2 x_1 \]
Solving for $x_1$, we find that 
\[ x_1 = \frac{d_1 + 2 d_2}{d_1 + 4 d_2} \implies x_2 + x_1 = \frac{3 d_1 + 6 d_2}{d_1 + 4 d_2} > 1 = x_3 \]
so $\bm{x}$ is not a feasible point. Thus we cannot have $z_1 + 2 z_2 = 0$, hence the complementary slackness conditions imply that $x_1 + x_2 = x_3$ .

Assume that $2 x_1 = x_2$. We claim this is only possible if and only if $n \leq 4$. We can write
\[ g_{\bm{d}}(\bm{x}) = \frac{1}{2} \sum_{i=1}^3 d_i(i x_1 - 1)^2 \]
and optimizing this as a function of $x_1$ only shows that
\[ x_1 = \frac{d_1 + 2 d_2 + 3 d_3}{d_1 + 4 d_2 + 9 d_3} \]
and this gives us the value of $x_2$ as well. We also have the condition that $0 \leq \lambda_1 = \frac{1}{3} (z_1 - z_2)$, so we must have $z_2 \leq z_1$. Writing $x_1, z_1$ and $z_2$ in terms of $n$ and performing a straightforward but tedious calculation, we find that $z_2 \leq z_1$ holds if and only if
\[ x_1 \geq \frac{d_2 - d_1}{2 d_2 - d_1} \iff \frac{-(n-4)(n+1)}{2n(3n+4)} \geq 0 \]
This shows the claim. 

When $2 x_1 = x_2$, we see that $\bm{x} = x_1 (1, 2, 3)^{\Transpose}$. Since the weights of the canonical filtration are defined up to scaling, this gives the claim in the $n \leq 4$ case. When $n > 4$, and hence $x_2 > 2 x_1$, the complementary slackness conditions implies that $z_1 = z_2$. This forces $z_1 = z_2 = \lambda_2$ so, recalling that $z_3 = -\lambda_2$, we solve for the $x_i$ in terms of $\lambda_2$ and find that
\[ x_1 = 1 - \frac{\lambda_2}{d_1}, \quad x_2 = 1 - \frac{\lambda_2}{d_2}, \quad x_3 = 1 + \frac{\lambda_2}{d_3} \]
Using the relation $x_1 + x_2 = x_3$, we solve for $\lambda_2$ and find that
\[ \lambda_2 = \frac{d_1 d_2 d_3}{d_1 d_2 + d_1 d_3 + d_2 d_3} \]
Substituting this into the expressions for the $x_i$ and writing the $d_i$ in terms of $n$, one finds that the $x_i$ are non-zero scalar multiples of the weights claimed in the statement.
\end{proof}

The example above shows that in general the weights of the canonical filtration depend on the dimensions of the filtered pieces in a non-obvious way. Also, for $n > 4$ we do not have $w(\mathfrak{m}^2) = 2 w(\mathfrak{m})$, so the canonical filtrations need not be homogeneous.

\subsection{Upper Triangular Matrix Algebras}

In this section we consider subalgebras of $\Mat_n(\basefield)$ which are block upper triangular. Let $\bm{n} = (n_1, \ldots, n_s) \in \bN^s$ such that $\sum_{i=1}^s n_i = n$. We denote by $A_{\bm{n}}$ the algebra whose elements are matrices of the form
\[ 
\begin{bmatrix}
A_{1,1} & A_{1,2} & \ldots & A_{1, s} \\
0 & A_{2,2} & \ldots & A_{2, s} \\
\vdots & \vdots & \ddots & \vdots \\
0 & 0 & \ldots & A_{s, s} \\
\end{bmatrix}
\]
where $A_{i,j} \in \Mat_{n_i, n_j}(\basefield)$. We shall call the entries in $A_{i,j}$ elements of the the $(i,j)$-th block. If $i \neq j$, we say that elements of $A_{i,j}$ are off-diagonal, while elements of $A_{i,i}$ are diagonal. Note that there is a group homomorphism $(\basefield^\ast)^n \to \Aut(A_{\bm{n}})$ given by the assignment $(t_1, \ldots, t_n) \mapsto \diag(t_1, \ldots, t_n)$ and thus we can apply our strategy.

\begin{prop} \label{prop: adapted basis for upper triangular matrix algebras}
The basis of standard matrix units $E_{i,j} \in A_{\bm{n}}$ forms an adapted basis of the canonical filtration. Furthermore, elements in the same off-diagonal blocks all have the same weight and the elements of the diagonal blocks all have weight $0$.
\end{prop}

\begin{proof}

Note that $A_{\bm{n}}$ is $\bZ^n$ graded. Indeed, for $1 \leq i,j \leq n$ if we set $\deg(E_{i,j}) = \bm{e}_j - \bm{e_i}$ one readily sees that this defines a grading. In particular, this shows that the matrix units $E_{i,j}$ with $j > i$ are part of an adapted basis for the canonical filtration.

Next we show that all elements of the same block have the same weight in the canonical filtration. We first consider the off-diagonal block in position $(i,j)$ where $j > i$. To ease notation, we now let $e_{l,q}$ in $\Mat_{n_i, n_j}(\basefield)$ denote the matrix unit relative to the $(i,j)$-th block. That is, if we set
\[ N_i := \sum_{a = 0}^{i-1} n_a, \quad N_j := \sum_{b=0}^{j-1} n_b \]
then we have $e_{l,q} = E_{N_i + s, N_j + t}$ in the ambient matrix algebra $A_{\bm{n}}$. Since $e_{l,q}$ is in the $(i,j)$-th block, multiplying on the right (respectively left) by a block diagonal matrix $\diag(A_1, \ldots, A_s)$ will give the matrix $e_{l,q} A_j$ (respectively $A_i e_{l,q}$) in the $(i,j)$-th block. Let $\pi_{st} \in \Mat_{n_j}(\basefield)$ be the permutation matrix corresponding to the transposition $(st)$, and set $P_{st, m} = \diag(I, \ldots, \pi_{st}, \ldots, I)$, where the $I$'s denote identity matrices of the appropriate size and $\pi_{st}$ occurs at position $m$. Observe that $\varphi_{st, j}(x) = P_{st, m} x P_{st, m}$ is an automorphism of $A_{\bm{n}}$, and hence preserves the canonical filtration by proposition (\ref{prop: action of automorphisms and derivations on canonical filtrations}). Now we compute $\varphi_{1 q, j}(e_{l,q})$. Since the block diagonal matrix $P_{1 q, j}$ has the identity matrix in the $i$-th block-diagonal spot, and $\pi_{1 q}$ in the $j$-th block diagonal spot, and $e_{l,q} \pi_{q 1} = e_{l,1}$, we find that
\[ \varphi_{1 q, j}(e_{l,q}) = P_{1 q, j} e_{l,q} P_{1 q, j} = P_{1 q, j} e_{l,1} = e_{l,1} \]
A similar computation shows that $\varphi_{1 l, i}(e_{l,1}) = e_{1, 1}$. This shows that the matrix unit $e_{l,q}$ in the $(i,j)$-th block is in $F(A)^{\Can}_m$ for some $m$ if and only if $e_{1,1}$ is in $F(A)^{\Can}_m$. Thus every matrix unit in the $(i,j)$-th block must have the same weight.

Now we consider the $(i,i)$-th diagonal block. If $e_{l,q}$ is a matrix unit in this block, then 
\[ \varphi_{lq, i}(e_{l,q}) = P_{lq,i} e_{l,q} P_{lq,i} = P_{lq,i} e_{l,l} = e_{q,l} \] 
since left and right multiplication by $P_{lk,i}$ effects rows and columns for diagonal blocks. This shows that $e_{l,k} \in F^{\Can}_m$ if and only if $e_{k,l} \in F^{\Can}_m$. Therefore $e_{l,k} e_{l,k} = e_{l,l} \in F^{\Can}_{2m}$, which implies that $m \leq 0$ since $e_{l,l}$ is not nilpotent. Thus any element in a diagonal block must have non-positive weight.

We now show that all elements in a diagonal block must have weight $0$ in the canonical filtration. Let $w$ be the weight function corresponding to the canonical filtration of $A_{\bm{n}}$. Define a new function $v$ by setting $v(E_{i,j}) = w(E_{i,j})$ if $E_{i,j}$ is an element of an off-diagonal block and $v(E_{i,j}) = 0$ if $E_{i,j}$ is an element of a diagonal block. We claim that this defines a weighted algebra basis for $A_{\bm{n}}$. Consider the product $E_{i,j} E_{l,q}$. We wish to show that
\[ v(E_{i,j}) + v(E_{l,q}) \leq v(E_{i,j} E_{l,q}) \]
If $l \neq j$ then $E_{i,j} E_{l, q} = 0$, so $v(E_{i,j} E_{l,q}) = \infty$ and there is nothing to check, so we can assume that $l = j$. If both $E_{i,j}$ and $E_{j,q}$ are elements of off-diagonal blocks then $E_{i,j} E_{j,q} = E_{i,q}$ is also an element of an off-diagonal block, so we have
\[ v(E_{i,j}) + v(E_{j,q}) = w(E_{i,j}) + w(E_{j,q}) \leq w(E_{i,q}) = v(E_{i,q}) \]
If $E_{i,j}$ is an element of a diagonal block while $E_{j,q}$ is an element of an off-diagonal block, then $E_{j,q}$ and $E_{i,q}$ are elements of the same off-diagonal block, hence $w(E_{j,q}) = w(E_{i,q})$. Thus
\[ v(E_{i,j}) + v(E_{j,q}) = 0 + w(E_{j,q}) = w(E_{i,q}) = v(E_{i,q}) \]
A similar argument shows that $v(E_{i,j}) + v(E_{j,q}) \leq v(E_{i,q})$ if $E_{j,q}$ belongs to a diagonal block while $E_{i,j}$ belongs to an off-diagonal block. In the last case, if both $E_{i,j}$ and $E_{j,q}$ belong to a diagonal block, then they necessarily belong to the same diagonal block. Hence their product $E_{i,q}$ belongs to the same diagonal block. Thus we have $v(E_{i,j}) = v(E_{j,q}) = v(E_{i,q}) = 0$, so the weight inequality is trivially satisfied in this case. This shows the claim that the weight $v(E_{i,j})$ define a weighted algebra basis.

Now we show that the weights $v(E_{i,j})$ define the weights of the canonical filtration. Let $G_\bullet$ denote the filtration corresponding to the weights $v(E_{i,j})$, and let $w(E_{i,j})$ denote the weights corresponding to the canonical filtration. The fact that $v(E_{i,j}) = 0 \geq w(E_{i,j})$ for all matrix units $E_{i,j}$ in diagonal blocks yields the inequalities
\[ \sum_{i,j} v(E_{i,j}) \geq \sum_{i,j} w(E_{i,j}), \quad \sum_{i,j} v(E_{i,j})^2 \leq \sum_{i,j} w(E_{i,j})^2 \]
It follows that
\[ \nu(A_{\bm{n}}, G_\bullet) = \frac{- \sum_{i,j} v(E_{i,j})}{\sum_{i,j} v(E_{i,j})^2} \leq \frac{- \sum_{i,j} w(E_{i,j})}{\sum_{i,j} w(E_{i,j})^2} = \nu(A_{\bm{n}}, F(A_{\bm{n}})^{\Can}_\bullet) \]
This shows that $F(A_{\bm{n}})^{\Can}_\bullet = G_\bullet$, and in particular the matrix units $E_{i,j}$ in the diagonal blocks must have weight $0$.
\end{proof}

Using proposition (\ref{prop: adapted basis for upper triangular matrix algebras}), the optimization problem computing the canonical filtration is as follows. For $1 \leq i \leq j \leq s$, let $d_{i,j} = n_i n_j$, which is the dimension of the $(i,j)$-th block of $A_{\bm{n}}$, and let $w_{i,j}$ denote the variable which is the weight of elements in the $(i,j)$-th block. We want to minimize
\[ f(\bm{w}) = \frac{- \sum_{i,j}d_{i,j} w_{i,j}}{\sum_{i,j} d_{i,j} w_{i,j}^2} \]
subject to the constraint $C \bm{w} \leq \bm{0}$ where $C$ is the matrix whose rows are of the form $\bm{e}_i + \bm{e}_j - \bm{e}_{i + j}$ for $1 \leq i,j \leq s$ and $i + j \leq s$.

We can explicitly solve this for $A_{\bm{1}}$ i.e. the algebra of upper triangular matrices in $\Mat_n(\basefield)$.

\begin{example} \label{example: canonical filtration of upper triangular matrices}
The canonical filtration of $A_{\bm{1}}$, the algebra of upper triangular matrices, is induced by the grading where the matrix units are the homogeneous elements and $\deg(E_{ij}) = j - i$.
\end{example}

\begin{proof}
We first compute the normalized weights. Set
\[ L := \sum_{1 \leq i,j \leq n} (j - i) = \sum_{l = 1}^{n-1} l (n-l) = \frac{n(n^2 - 1)}{6}, \quad B:= \sum_{l = 1}^{n-1} l^2 (n-l) = \frac{n^2(n^2 - 1)}{12} \]
The normalized weight of $E_{i,j}$ is thus $w(E_{i,j}) = \frac{L}{B} (j - i) = \frac{2(j - i)}{n}$. To verify these are the maximizing weights, we need only verify the condition of Farkas' lemma since these weights define a grading. We need to show that for any set of real numbers $y_{i,j}$, $1 \leq i \leq j \leq k \leq n$, such that $y_{i,j} + y_{j,k} \leq y_{i,k}$ (coming from the constraints on the weights) we have
\[ \sum_{i \leq j} (1 - w(E_{i,j})) y_{i,j} = \frac{1}{n} \sum_{i \leq j} (n - 2(j - i)) y_{i,j} \leq 0 \]
Note that for terms where $i = j$ the term is $n y_{i,i} \leq 0$, and we know that $y_{i,i} \leq 0$, so it suffices to consider the above sum but with terms where $i < j$. We claim that
\begin{equation} \label{eqn: upper triangular sum farkas lemma}
\sum_{i < j} (n - 2(j - i)) y_{i,j} = \sum_{i < j < q} y_{i,j} + y_{j,q} - y_{i,q} 
\end{equation}
which, in light of the inequalities $y_{i,j} + y_{j,q} \leq y_{i,q}$, will conclude the proof. We will show that for any pair $(i,j)$ the coefficient of $y_{i,j}$ is the same in both sums in (\ref{eqn: upper triangular sum farkas lemma}). Note that $y_{i,j}$ will appear in the terms of the right-hand sum with a $+1$ coefficient for every product of the form $E_{i,j} E_{j,q}$ or $E_{q,i} E_{i,j}$, and will appear with a $-1$ coefficient for every product of the form $E_{i,q} E_{q,j}$. The number of products of the former type is $(n - j) + (i - 1)$, while the number of products of the later type is given by $j - i - 1$, so the total coefficient of $y_{i,j}$ in the right-hand sum is $(n - j) + (i - 1) - (j - i - 1) = n - 2(j - i)$. This is exactly the coefficient of $y_{i,j}$ in the left-hand sum, as claimed. Thus these weights define the canonical filtration.
\end{proof}

\subsection{Lie Algebras}

In this section we compute the canonical filtrations of Lie algebras which are amenable to our methods.

\subsubsection{Filiform Lie Algebras}
A \emph{filiform} Lie algebra $\mathfrak{g}$ of dimension $n + 1$ whose lower central series defined as $\mathfrak{g}_{m}$ satisfies $\dim(\mathfrak{g}_{m}) = n - m - 1$ for $ 0 \leq m \leq n - 1$. By \cite[Proposition 2]{remm2017filiform}, every filiform lie algebra has a basis of the following form.

\begin{defn} \label{defn: Vergne basis}
A Vergne basis of a filiform Lie algebra is a basis $\{ x_0, x_1, \ldots, x_n \}$ such that
\begin{enumerate}
    \item[(i)] $[x_0, x_i] = x_{i + 1}$ for $i = 1 \leq i \leq n -1$.
    \item[(ii)] $[x_1, x_{n-1}] = 0$
    \item[(iii)] $[x_i, x_j] \subseteq \Span(x_{i+j}, \ldots, x_n)$.
\end{enumerate}
\end{defn}

It follows from this definition that $\mathfrak{g}_i = \Span(x_i, \ldots, x_n)$ for $1 \leq i \leq n - 2$. These bases allow us to give some structural results about the canonical filtration of filiform Lie algebras.

\begin{lemma} \label{lemma: adapted basis for filiform lie algebras}
Let $\mathfrak{g}$ be a filiform Lie algebra with Vergne basis $x_0, \ldots, x_n$. Suppose that for some $m$ we have
\[  \sum_{j = 1}^r c_{i_j} x_{i_j} \in F(\mathfrak{g})^{\Can}_m  \]
where $1 \leq i_1 < \ldots < i_r \leq n$ and $c_i \neq 0$ for all $i \in S$. Then $x_{i_1}, \ldots, x_{i_r} \in F(\mathfrak{g})^{\Can}_m$.
\end{lemma}

\begin{proof}
Let $y = \sum_{i \in S} c_i x_i$. Since the canonical filtration is invariant under derivations by proposition (\ref{prop: action of automorphisms and derivations on canonical filtrations}), we have
\[ \ad_{x_0}^s(y) = \sum_{j = 1}^r c_{i_j} x_{i_j + s} \in F(\mathfrak{g})^{\Can}_m  \]
where we set $x_{t} = 0$ if $t > n$. Let $C$ be the square matrix with $n - i_1 + 1$ rows and columns such that the entry $C_{pq}$ is the coordinate of $x_{q - 1 + i_1}$ in the basis expansion of $\ad_{x_0}^{p- 1}(y)$. It follows from the calculation above that $C_{pq} = 0$ if $p < q$ and that $C_{pp} = c_{i_1} \neq 0$ for all $p$, hence $C$ is invertible. Now let $\bm{y}$ and $\bm{x}$ be the length $n - i_1 + 1$ column vectors with entries in $\mathfrak{g}$ given by
\[ \bm{x}_{q} = x_{i_1 + q - 1}, \quad \bm{y}_{p} = \ad_{x_0}^{p - 1}(y) \]
Then $\bm{y} = C \bm{x}$, and since $C$ is invertible we have $\bm{x} = C^{-1} \bm{y}$. Thus each $x_{i_1}, \ldots, x_n$ can be expressed as a $\basefield$-linear combination of the elements $\ad_{x_0}^{p}(y)$, which implies that $x_{i_1}, \ldots, x_n \in F(\mathfrak{g})^{\Can}_m$.
\end{proof}

\begin{prop} \label{prop: canonical filtration of filiform lie algebras}
The canonical filtration of a filiform Lie algebra $\mathfrak{g}$ is a refinement of its lower central series.
\end{prop} 

\begin{proof}
Let $\mathfrak{g}_i$ denote the lower central series of $\mathfrak{g}$, with $\mathfrak{g}_0 = \mathfrak{g}$. If $x_0, x_1, \ldots, x_n$ is a Vergne basis of $\mathfrak{g}$ then for $i \geq 1$ we have $\mathfrak{g}_i = \Span(x_{i+1}, \ldots, x_n)$. Let $i_m$ denote the smallest integer $j$ such that there exists $y \in F(\mathfrak{g})^{\Can}_m$ such that $y = \sum_{t \geq j} c_t x_t$ with $c_j \neq 0$. If $i_m \geq 1$, we claim that $\Span(x_{i_m}, \ldots, x_n) = F(\mathfrak{g})^{\Can}_m$. By lemma (\ref{lemma: adapted basis for filiform lie algebras}) we have $x_{i_m}, \ldots, x_n \in F(\mathfrak{g})^{\Can}_m$ and hence $\Span(x_{i_m}, \ldots, x_n) \subseteq F(\mathfrak{g})^{\Can}_m$. If $y = \sum_{t \geq j} c_t x_t \in F(\mathfrak{g})^{\Can}_m$, then by the minimality of $i_m$ we have $j \geq i_m \geq 1$, so $F(\mathfrak{g})^{\Can}_m \subseteq \Span(x_{i_m}, \ldots, x_n)$, proving the claim. This shows that if $i_m \geq 2$ then $F(\mathfrak{g})^{\Can}_m = \mathfrak{g}_{i_m - 1}$, hence the canonical filtration refines the lower central series.
\end{proof}

For $n \geq 1$, we define the model filiform Lie algebra $M_n$ of dimension $n + 1$ as the filiform Lie algebra such that the Lie brackets in (i) of definition (\ref{defn: Vergne basis}) are the only non-zero brackets; see also \cite[\S 1]{remm2017filiform}. Note that $M_n$ is a $\bZ$-graded Lie algebra with grading $\deg(x_i) = i$ for $i \geq 1$ and $\deg(x_0) = 1$. This doesn't immediately fit into the strategy outlined in \autoref{subsection: Strategy}, as the degree one space has dimension two. However, the next lemma shows that $M_n$ has enough automorphisms to ``separate'' the degree one elements of a Vergne basis, so we can apply our strategy.

\begin{lemma} \label{lemma: an automorphism of M_n}
For any $c \in \basefield$, let $\varphi_c \in \End(M_n)$ be given by
\[ \varphi_c(x_0) = x_0 + c x_1, \quad \varphi_c(x_i) = x_i, \quad i > 0 \]
on a Vergne basis $\{ x_0, x_1, \ldots, x_n \}$. Then $\varphi_c$ is an automorphism of $M_n$. Consequently, a Vergne basis is adapted to the canonical filtration.
\end{lemma}

\begin{proof}
Observe that if $i, j >0$ from the definition we have that
\[ \varphi_c([x_i, x_j]) = 0 = [\varphi_c(x_i), \varphi_c(x_j)] \]
\[ \varphi_c([x_0, x_i]) = x_{i+1} = [x_0 + c x_1, x_i] = [\varphi_c(x_0), \varphi(x_i)] \]
Since $\varphi_c \circ \varphi_{-c} = \id$, $\varphi_c$ is an automorphism. Given the grading on $M_n$ above and corollary (\ref{coroll: canonical filtration is a direct sum of graded pieces}), to show that a Vergne basis is adapted to the canonical filtration it is enough to show that if $a x_0 + b x_1 \in F^{\Can}_m(M_n)$ then $a x_0, b x_1 \in F^{\Can}_m(M_n)$. If $a = 0$ there is nothing to show, so after scaling we can assume that $a = 1$ and $x_0 + b x_1 \in F^{\Can}_m(M_n)$. But then $\varphi_{-b}(x_0 + b x_1) = x_0 \in F^{\Can}_m(M_n)$ and hence $(x_0 + b x_1) - x_0 = b x_1 \in F^{\Can}_m(M_n)$ as required.
\end{proof}

\begin{example} \label{example: model filiform Lie algebras}
The canonical filtration of $M_n$ is given by the weights
\[ w(x_0) = 12, \quad w(x_i) = n^3 - 7n + 18 + 12(i-1) \quad 1 \leq i \leq n \]
\end{example}

\begin{proof}
For $0 \leq i \leq n$, let $\bm{e}_i$ be the standard basis vectors of $\bR^{n+1}$ and let $\bm{w}$ be the vector whose $i$-th entry is $w(x_i)$. The constraint matrix $C$ is the $n -1 \times (n + 1)$ matrix whose $i$-th row is $C_i = \bm{e}_0^{\Transpose} + \bm{e}_i^{\Transpose} - \bm{e}_{i+1}^{\Transpose}$. If we set $z_i = 1 - w_i$, the first KKT condition says that
\[ z_0 = \sum_{i=1}^{n-1} \lambda_i, \quad z_1 = \lambda_1, \quad z_n = - \lambda_{n-1}, \quad z_{i + 1} = -\lambda_i + \lambda_{i+ 1} \quad 1 \leq i \leq n - 2 \]
Solving for the $\lambda_i$ in terms of the $z_i$, we find that for $1 \leq i \leq n-1$,
\[ \lambda_i = \sum_{j=1}^i z_j = i - \sum_{j=1}^i w_j \]
Now we show that if all constraints hold with equality, so that $w_{i + 1} = w_i + w_0$, then we can solve for $w_0$ and $w_1$  and show that all KKT conditions are satisfied. A quick computation shows that in this case we have $w_i = w_1 + (i-1) w_0$ for $1 \leq i \leq n$. Taking $i = n - 1$ above, we find that
\[ \lambda_{n-1} = - z_n = \sum_{j=1}^{n-1} z_j \implies \sum_{j=1}^n z_j = 0 \iff n = \sum_{j = 1}^n w_j \]
Rewriting this in terms of $w_0$ and $w_1$ only, we have
\begin{equation} \label{eqn: filiform equality 1}
n = \sum_{j=1}^n w_1 + (j-1)w_0 = n w_1 + \frac{n(n-1)}{2} w_0
\end{equation}
From the expression above for $z_0$, we find
\[ 1 - w_0 = \sum_{i=1}^{n-1} \left( i - \sum_{j=1}^i w_j \right) = \frac{n(n-1)}{2} - \sum_{i=1}^{n-1} \sum_{j=1}^i w_1 + (j-1) w_0 \]
\[ = \frac{n(n-1)}{2} - \frac{n(n-1)}{2} w_1 - \frac{n(n-1)(n-2)}{6} w_0 \]
Solving equation (\ref{eqn: filiform equality 1}) for $w_1$ and substituting this into the equation above, we find that
\[ w_0 = \frac{12}{n^3 - n + 12} \]
and hence
\[ w_1 = 1 - \frac{n-1}{2} w_0 = \frac{n^3 - 7n + 18}{n^3 - n + 12} \]
Clearing denominators then gives the weights in the statement. Lastly, we check that $\lambda_i \geq 0$ for $1 \leq i \leq n - 1$:
\[ \lambda_i = i - \sum_{j=1}^i ( w_1 + (j-1)w_0 ) = i - \sum_{j=1}^i \frac{n^3 - 7n + 18 + 12(j-1)}{n^3 - n + 12} \]
\[ = i \left(1 - \frac{n^3 - 7n + 18 + 6(i-1)}{n^3 - n + 12} \right) \]
The expression in the parentheses is smallest when $i = n - 1$, and one can check it is greater than zero in this case. This shows that $\lambda_i \geq 0$.
\end{proof}

\subsubsection{Borel Subalgebras of Semisimple Lie Algebras}
In this subsection we assume that $\basefield$ is algebraically closed and of characteristic $0$ so that we have root space decompositions of semisimple Lie algebras. Thus for a semisimple Lie algebra $\mathfrak{g}$, we can choose a Cartan subalgebra $\mathfrak{h}$ and a corresponding set of roots $\Phi$ such that 
\[ \mathfrak{g} = \bigoplus_{\alpha \in \Phi} \mathfrak{g}_\alpha \]
If $\Phi^+$ is a set of positive roots, we form the associated nilpotent subalgebra
\[ \mathfrak{n} = \bigoplus_{\alpha \in \Phi^+} \mathfrak{g}_\alpha \]
and Borel subalgebra $\mathfrak{b} = \mathfrak{h} \oplus \mathfrak{n}$. Since the root system $\Phi$ defines a $\bZ^r$-grading on the algebras $\mathfrak{b}$ and $\mathfrak{n}$ and $\dim(\mathfrak{g}_\alpha) = 1$ for each root $\alpha$, our strategy for computing canonical filtrations applies. Furthermore, by theorem (\ref{thm: sufficient to compute canonical filtration in the nilpotent case}) it suffices to compute the canonical filtration of $\mathfrak{n}$.
 
 We will illustrate by computing the canonical filtration of a Borel subalgebra of $\mathfrak{sl}_{n+1}$. One choice of simple roots in $\bR^{n+1}$ is the set $\Phi$ of vectors of the form $\alpha_i = \bm{e}_{i} - \bm{e}_{i + 1}$. The positive roots $\Phi^+$ are then vectors of the form $\beta_{i,j} := \sum_{l=i}^{j-1} \alpha_l$ for $1 \leq i < j \leq n + 1$. Note that this defines a $\bZ^{n+1}$-grading on $\mathfrak{sl}_n$.

\begin{prop} \label{prop: canonical filtration of borels of sl_n+1}
Let $\mathfrak{n} \subseteq \mathfrak{sl}_{n+1}$ be as above, and let $\Phi_m \subseteq \Phi^+$ be the set of roots $\beta_{i, j}$ with $j - i \geq m$. Then
\[ F(\mathfrak{n})^{\Can}_m = \bigoplus_{\alpha \in \Phi_m} (\mathfrak{sl}_{n+1})_\alpha \]
\end{prop}

\begin{proof}
Let $G_\bullet$ be the filtration on the right-hand side of the equality in the proposition statement. First we compute the normalized weights of $G_\bullet$. Set
\[ L = \sum_{m=1}^n m \dim(G_m / G_{m+1}) = \sum_{m = 1}^n m(n - m + 1) =  \frac{n(n+1)(n+2)}{6} \]
\[ B = \sum_{m=1}^n m^2 \dim(G_m / G_{m+1}) = \sum_{m = 1}^n m^2(n - m + 1) =  \frac{n(n+1)^2(n+2)}{12} \]
The basis vector $E_{i,j}$ of the root space $\beta_{i, j}$ has weight $j - i$ in the filtration $G_\bullet$, so it normalized weight is
\[ \frac{L(j-i)}{B} = \frac{2(j - i)}{n + 1} \]
To confirm that these weights maximize the numerical invariant, it suffices to check the condition of Farkas' lemma (\ref{lemma: Farkas' lemma}). Let $y_{i,j}$ be a sequence of real numbers satisfying $y_{i,j} + y_{j,q} \leq y_{i,q}$ for $1 \leq i < j < q \leq n + 1$. We must verify that
\[ \sum_{1 \leq i < j \leq n + 1} \left(1 - \frac{2(j - i)}{n + 1} \right) y_{i,j} = \frac{1}{n+1} \sum_{1 \leq i < j \leq n + 1} \left(n + 1 - 2(j - i)\right) y_{i,j} \leq 0 \]
We claim that
\[
\sum_{1 \leq i < j \leq n + 1} \left(n + 1 - 2(j - i)\right) y_{i,j} = \sum_{1 \leq i < j < q \leq n + 1} y_{i,j} + y_{j,q} - y_{i,q}
\]
which, in light of the inequalities $y_{i,j} + y_{j,q} \leq y_{i,q}$, is enough to verify the inequality above. The proof of the claim is similar to that of examples (\ref{example: jordan algebras}) and (\ref{example: canonical filtration of upper triangular matrices}). For each pair $(i,j)$, the term $y_{i,j}$ appears with a coefficient of $+1$ in the right-hand sum a total of $(n + 1 - j) + (i - 1)$ times, and with a coefficient of $-1$ a total of $j - i - 1$ times. Thus the coefficient of $y_{i,j}$ in the right-hand sum is $n + 1 - 2(j - i)$, which is exactly the coefficient of $y_{i,j}$ in the left-hand sum. This verifies the conditions of Farkas' lemma and proves that $G_\bullet$ is the canonical filtration as claimed.
\end{proof}

\section{The Structure of Graded-Semistable Algebras} \label{section: characterizations of graded-semistable algebras}

In general, describing the canonical filtration of an arbitrary algebra is a difficult problem. However, for the associated graded algebra, $\Gr(F(A)^{\Can}_\bullet)$ we can say a bit more. We recall the notion of graded-semistability in our set-up, where we consider the stability of a graded algebra $A_\bullet$ with respect to a modified stability condition. We recover some of the results from \cite{zhang2023moment} on the algebraic structure of the critical points of the norm squared of the moment map, which corresponds in our setting to graded-semistable algebras. However, we do not recover their results on the compatibility of the graded component with an inner product on the underlying vector space of the algebra, since $\basefield$ may have positive characteristic where defining inner products requires some care.

We define the \emph{split filtration} $F(A_\bullet)^{\Split}_\bullet$ of a graded algebra $A_\bullet$ is the filtration induced by the grading, i.e.
\[ F(A_\bullet)^{\Split}_m = \bigoplus_{d \geq m} A_d \]

\begin{defn} \label{defn: graded-semistability}
A graded algebra $A_\bullet$ is graded-semistable if the split filtration is the canonical filtration.
\end{defn}

There is a geometric invariant theory interpretation of graded-semistability, which we will recall from \cite[\S 2.3]{hoskins2014stratifications}. Note that \cite{hoskins2014stratifications} considers the case of GIT on an affine space, but this generalizes to affine schemes with minimal modifications. We say that a filtration $F(A_\bullet)_\bullet$ is \emph{compatible} with the grading of $A_\bullet$ if it is given by a bounded $\bZ$-filtration $F(A_n)_\bullet$ for each graded component $A_n$ and satisfies
\[ F(A_{n})_p \cdot F(A_{n'})_q \subseteq F(A_{n + n'})_{p + q}\]
for all $n, n', p, q \in \bZ$. Let $X$ be either $\ModAlg$ or $\ModLieAlg$ and let $\lambda$ be a non-trivial one-parameter subgroup of $G := \GL(V)$ which defines the canonical filtration of some set of points in $X$. Set
\[ Z_{\lambda} = \{ \mu_A \in X \mid \lambda \text{ defines a grading on } A \text{ and } F(\lambda)_\bullet = F(A_\bullet)^{\Can}_\bullet \}  \]
\[ X^\lambda = \{ \mu_A \in X \mid \lambda \text{ defines a grading on } A \} \]
\[ G_\lambda = \{ g \in G \mid g \text{ commutes with } \lambda \} \]
Thus $Z_\lambda$ is the set of graded-semistable algebras with grading defined by $\lambda$. Observe that $G_\lambda$ acts on $X^\lambda$. We can linearize the trivial bundle with respect to this action via the character $\rho_{\lambda}$ of $G_\lambda$ defined to be
\[ \rho_\lambda(\lambda') := || \lambda ||^2 \det(\lambda')^{-1} - \det(\lambda)^{-1} (\lambda, \lambda') \]
where $\lambda'$ is a one-parameter subgroup of $G_\lambda$ and $(\dash, \dash)$ denotes the pairing of one-parameter subgroups. Since $\det(\lambda)^{-1}$ and $|| \lambda ||^2$ are positive, this does in fact define a character. Then by \cite[Proposition 2.18]{hoskins2014stratifications}, the locus $(X^\lambda)^{\rho_\lambda - \semistable}$ of points in $X^\lambda$ which are semistable with respect to $\rho_\lambda$ is precisely $Z_\lambda$. In other words, $A$ is graded-semistable if and only if it is $\rho_\lambda$-semistable. Furthermore, thanks to \cite[Theorem 5.4.4]{halpern2014structure}, graded-semistability gives us a criterion to determine if an arbitrary filtration is the canonical filtration:

\begin{thm}[The Recognition Theorem] \label{thm: the recongnition theorem}
Let $F_\bullet$ a destabilizing filtration of an algebra $A$. Then $F_\bullet$ is the canonical filtration if and only if the associated graded algebra $\Gr(F_\bullet)$ of $F_\bullet$ is graded-semistable.
\end{thm}

We now show that the canonical filtration of a graded-semistable algebra is, in the appropriate sense, dual to a certain linear functional. First we need some set up. Recall that to give a $\bZ$-grading on an algebra $A$ is the same as giving a $\bG_{m}$ action on $A$, which we can identify with an element $\varphi \in \Aut(A)$ (also identified with a one-parameter subgroup). We will call such an automorphism a \emph{grading operator}. After choosing the appropriate basis, we can identify a grading operator $\varphi$ with a diagonal matrix $\diag(t^{m_1}, \ldots, t^{m_d})$ where $m_i \in \bZ$. Thus the $m_i$ correspond to the degrees of the basis elements. We call a grading operator trivial if the grading is concentrated in degree $0$, or equivalently $\varphi = \id$. We define the \emph{grading trace} of $\varphi$ as
\[ \GTrace(\varphi) := \sum_{j=1}^d m_j \]

Now fix a grading operator $\varphi$, and let $C(\varphi)$ denote the set of all grading operators $\psi$ which commute with $\varphi$. Thus $C(\varphi)$ is the set of gradings which refine the grading induced by $\varphi$. We define an inner product $B$ on $C(\varphi)$ by
\[ B(\psi_1, \psi_2) := \GTrace(\psi_1 \circ \psi_2\textbf{}) = \sum_{j = 1}^d m_{1,j} m_{2,j} \]
where we have identified $\psi_i$ with $\diag(t^{m_{i,1}}, \ldots, t^{m_{i,d}})$. We denote by $|| \dash ||_B$ the norm induced by $B$. Note that $\psi \circ \varphi$ is also a grading operator for $\psi \in C(\varphi)$, since commuting diagonalizable linear maps are simultaneously diagonalizable. It follows from the definitions that
\[ \numinvariant(A_\bullet, F(\psi)_\bullet ) = - \frac{\GTrace(\psi)}{|| \psi ||_B} \]
where $F(\psi)_\bullet$ is the split filtration defined by the grading induced by $\psi$.

Observe that $C(\varphi)$ has the structure of a $\bZ$ module. Indeed, if we write $\psi \in C(\varphi)$ as a diagonal matrix $\diag(t^{m_1}, \ldots, t^{m_d})$ in the weight space basis for $\varphi$, then we can identify $\psi$ with the vector $(m_1, \ldots, m_d)$ and the module structure is then the obvious one. Thus, we can extend scalars to $\bR$ so that $C(\varphi)_{\bR} := C(\varphi) \otimes_{\bZ} \bR$ has the structure of an $\bR$ vector space. By abuse of notation, we will also let $\GTrace$ and $B$ denote their extensions to $C(\varphi)_{\bR}$. For $r \in \bR$ and $\psi(t) \in C(\varphi)$, which we identify with $\diag(t^{m_1}, \ldots ,t^{m_d})$ with $m_i \in \bZ$, the operator $\psi(t^r)$ defines a $\bR$-grading on $A$ such that the elements of degree $r m_i$ are is $A_{m_i}$.

The following proposition says that, up to scaling, the grading operator of a graded-semistable algebra $\varphi_{\Can}$ is dual to the inner product $B$ on $C(\varphi_{\Can})$. 

\begin{prop} \label{prop: riez element for commuting grading operators}
Let $A$ be a graded algebra which is graded-semistable, $\varphi$ be a non-trivial grading operator, and $\varphi_{\Can}$ the grading operator defining the canonical filtration. Then $\varphi_{\Can} \in C(\varphi)$ and
\[ \GTrace(\varphi_{\Can} \circ \psi) = \frac{|| \varphi_{\Can}||_B^2}{\GTrace(\varphi_{\Can})} \GTrace(\psi) \]
for any $\psi \in C(\varphi)$.
\end{prop}

\begin{proof}
Since $B$ is a non-degenerate bilinear form on $C(\varphi)_{\bR}$, there is some $\phi_{\riesz} \in C(\varphi)_{\bR}$ such that for any $\psi \in C(\varphi)_{\bR}$
\[ B(\phi_{\riesz}, \psi) = \GTrace(\psi) \]
Taking $\psi = \phi_{\riesz}$, we see that
\[ || \phi_{\riesz} ||_B^2 = \GTrace(\phi_{\riesz}) \]
From the Cauchy-Schwarz inequality, we have
\[ B(\phi_{\riesz}, \psi) = \GTrace(\psi) \leq || \phi_{\riesz} ||_B \cdot || \psi ||_B \]
which implies that $\numinvariant(A_\bullet, F(\psi)) \geq \numinvariant(A_\bullet, F(\phi_{\riesz}))$ for any non-trivial $\psi \in C(\varphi)_{\bR}$. Thus $\phi_{\riesz}$ is the minimizer of $\numinvariant(A_\bullet, \dash)$ on $C(\varphi)_{\bR}$. Proposition (\ref{coroll: canonical filtration is a direct sum of graded pieces}) implies that the canonical filtration of $A$ is compatible with the grading defined by $\varphi$, which implies that $\varphi_{\Can}$ commutes with $\varphi$, hence $\varphi_{\Can} \in C(\psi)$. Thus we have
\[\numinvariant(A_\bullet, F(\varphi_{\Can})) \leq \numinvariant(A_\bullet, F(\phi_{\riesz})) \leq \numinvariant(A_\bullet, F(\varphi_{\Can}))\] 
By the uniqueness of the canonical filtration, we must have $\varphi_{\Can}(t) = \phi_{\riesz}(t^K)$ for some $K \in \bR$. 
Observe that
\[ || \varphi_{\Can} ||_B^2 = K^2 B(\phi_{\riesz}, \phi_{\riesz}) = K^2 \GTrace(\phi_{\riesz}) = K \GTrace(\varphi_{\Can}) \]
Since $\varphi$ is non-trivial, $\GTrace( \phi_{\riesz} ) \neq 0$ and we find that
\[ K = \frac{|| \varphi_{\Can} ||_B^2}{\GTrace(\varphi_{\Can})} \]
This implies the claim.
\end{proof}

We note that we can interpret proposition (\ref{prop: riez element for commuting grading operators}) in terms of the weights and norms of filtrations. Recalling that $F(\psi)_\bullet$ denotes the filtration induced by a grading operator $\psi$, it follows from the definitions that
\[ \wt(A, F(\psi)_\bullet) = \GTrace(\psi), \quad || F(\psi)_\bullet || = || \psi ||_B \]
so the conclusion of the proposition can be restated as
\[ \wt(A, F(\varphi_{\Can} \circ \psi)_\bullet ) = \frac{|| F(\varphi_{\Can})_\bullet||^2}{\wt(A, F(\varphi_{\Can})_\bullet )} \wt(A, F(\psi)_\bullet ) \]

We note that proposition (\ref{prop: riez element for commuting grading operators}) says that the canonical grading operator of a graded-semistable algebra is the analogue of a Nikolayevsky derivation $\delta_{\Nik}$ of an algebra $A$ over $\bC$ as defined in \cite[\S 4.1]{zhang2023moment}. This derivation has the property that for any derivation $\delta$,
\begin{equation} \label{eqn: Nikolayevsky property}
\Trace(\delta \circ \delta_{\Nik}) = \Trace(\delta)
\end{equation}
and that $\delta_{\Nik}$ is diagonalizable with rational eigenvalues.
For a diagonalizable derivation $\delta$ with integral eigenvalues, one can define a grading on $A$ by taking the graded components to be the eigenspaces of $A$. If we let $F(\delta)_\bullet$ denote the split filtration induces by such a derivation, one sees that
\[ \wt(A, F(\delta)_\bullet) = \Trace(\delta), \quad || F(\delta)_\bullet || = \sqrt{\Trace(\delta^2)} \]
Since $\delta_{\Nik}$ has rational eigenvalues, $\Delta_{\Nik} = N \delta_{\Nik}$ will have integral eigenvalues for some sufficiently divisible $N \in \bZ^+$, hence it will define a $\bZ$-filtration of $A$. Then we will have
\[ \Trace( \Delta_{\Nik}^2 ) = N^2 \Trace(\delta_{\Nik}) = N \Trace(\Delta_{\Nik}) \]
and thus 
\[ N = \frac{\Trace( \Delta_{\Nik}^2 )}{\Trace(\Delta_{\Nik})} = \frac{ || F(\Delta_{\Nik})_\bullet ||^2}{\wt(A, F(\Delta_{\Nik})_\bullet)} \]
Furthermore, one can show that the equality (\ref{eqn: Nikolayevsky property}) implies that
\[ \numinvariant(A, F(\delta)_\bullet) = -\frac{\Trace(\delta)}{\sqrt{\Trace(\delta^2)}} \geq -\frac{\Trace(\Delta_{\Nik})}{\sqrt{\Trace(\Delta_{\Nik}^2)}} = \numinvariant(A, F(\Delta_{\Nik})_\bullet) \]
and hence if $A$ is graded-semistable then the canonical filtration is given by $F(\delta_{\Nik})_\bullet$.

Now we give some structural results of graded-semistable associative algebras. Recall that the annihilator of an algebra $A$ is the ideal
\[ \Ann(A) = \{ a \in A \mid ab = ba = 0 \text{ for all } b\in A \} \]
If $A$ is graded, then $\Ann(A)$ is a homogeneous ideal, and we denote its $n$-th homogeneous component by $\Ann(A)_n$.

\begin{prop} \label{prop: annihilator of graded-semistable ring}
Let $A_\bullet$ be a graded-semistable algebra with non-trivial grading, and let $\Ann(A_\bullet)$ denote the annihilator of $A_\bullet$. Then $\Ann(A_\bullet)_n = 0$ if $n < 0$.
\end{prop}

\begin{proof}
We define a filtration $F_\bullet$ which is compatible with the grading of $A_\bullet$ via a weight function $w$. Choose a homogeneous basis $\cA$ for $\Ann(A_\bullet)$ and for every element $x \in \cA$, set $w(x) = |\deg(x)|$. Now for each $n$, choose a basis $\cB_n$ of $A_n \setminus \Ann(A_\bullet)_n$ and set $w(y) = n$ for all $y \in \cB_n$. Let the basis $\cB$ be the union of the sets $\cB_n$ and $\cA$. Since our basis consists of homogeneous elements, $F_\bullet$ is compatible with the grading by construction. We now show that for any $x,y \in \cB$ such that $xy \neq 0$, we have
\[ w(x) + w(y) \leq w(xy) \]
Note that in this case neither $x$ nor $y$ is in $\Ann(A_\bullet)$. Since $xy$ is a homogeneous element, if we write the basis expansion
\[ xy = \sum_{b \in \cB \setminus \cA} c_b b + \sum_{a \in \cA} c_a a \]
then for each $a$ and $b$ appearing with $c_a \neq 0$ and $c_b \neq 0$, we must have that $\deg(a) = \deg(b) = \deg(xy)$. Denote this common degree by $n$. For these $a$ and $b$, we have $w(a) = |n| \geq n$ and $w(b) = n$, so the minimum of the weights over all $a$ and $b$ with non-zero coefficients is $n$. Since $x,y \notin \Ann(A_\bullet)$, we have 
\[ w(x) + w(y) = \deg(x) + \deg(y) = \deg(xy) = n \]
Thus, for each $a$ and $b$ with non-zero coefficients above, $w(a)$ and $w(b)$ are greater than or equal to $w(x) + w(y)$, so $w(xy) \geq w(x) + w(y)$ as required.

From the definition of $F_\bullet$, we have
\begin{align*}
\wt(A_\bullet, F(A_\bullet)^{\Split}_\bullet) &= \sum_n n \dim(A_n) \\
&= \sum_n n \dim(A_n / \Ann(A_\bullet)_n) + \sum_n n \dim(\Ann(A_\bullet)_n) \\
&\leq \sum_n n \dim(A_n / \Ann(A_\bullet)_n) + \sum_n |n| \dim(\Ann(A_\bullet)_n) \\
&= \wt(A_\bullet, F_\bullet)
\end{align*}
A similar calculation shows that 
\[ || F(A_\bullet)^{\Split}_\bullet || = || F_\bullet || \]
It follows that $\numinvariant(A_\bullet, F(A_\bullet)^{\Split}_\bullet) \geq \numinvariant(A_\bullet, F_\bullet)$, and in fact this is an equality since $A_\bullet$ is graded-semistable. But the only way this is possible is if $\Ann(A_\bullet)_n = 0$ for $n < 0$, proving the claim.
\end{proof}

Now we give some results about the algebraic structure of graded-semistable algebras, and in particular recover the results of 
\cite[Theorem 4.10]{zhang2023moment}.

\begin{prop} \label{prop: graded-semistability and the canonical filtration}
Let $A_\bullet$ be a graded-semistable associative algebra with non-trivial grading. Let $\JRad(A_\bullet)$ and $Z(A_\bullet)$ denote the Jacobson radical and center of $A_\bullet$ respectively, and set
\[ A_\bullet^{-} = \bigoplus_{i < 0} A_i, \quad A_\bullet^{+} = \bigoplus_{i > 0} A_i \]
Then 
\begin{enumerate}
    \item[(i)] $A_\bullet^{-} \subseteq \JRad(A_\bullet) \cap Z(A_\bullet)$
    \item[(ii)] $A_\bullet^{+} \subseteq \JRad(A_\bullet)$.
    \item[(iii)] $\Ann(A_\bullet)_n = 0$ for $n < 0$.
\end{enumerate}
In particular, this holds for $A_\bullet := \Gr(F(A)^{\Can}_\bullet)$ where $A$ is an unstable associative algebra.
\end{prop}

\begin{proof}
Note that (iii) follows immediately from the graded-semistability of $A_\bullet$ and lemma (\ref{prop: annihilator of graded-semistable ring}). Since $\Gr(F(A)^{\Can}_\bullet)$ is graded-semistable by the recognition theorem (\ref{thm: the recongnition theorem}), this will give the final statement once we have shown everything else.

Since $A_\bullet$ is graded-semistable, the canonical filtration of $A_\bullet$ is the split filtration. First we show that $A_\bullet^{-} \subseteq Z(A_\bullet)$. Let $x \in A_n$ and $y \in A_m$, with $n < 0$. Consider the derivation $\delta_x(a) = xa - ax$. Clearly $\delta_x(y) \in A_{m + n}$. Furthermore, proposition (\ref{prop: action of automorphisms and derivations on canonical filtrations}) shows that 
\[ \delta_x(y) \in A_{m +n} \bigcap \bigoplus_{i \geq m} A_i \]
Since $m + n < m$, this intersection is $0$. Thus $x y = y x$ and $A_\bullet^{-} \subseteq Z(A_\bullet)$. Now we show that $A_\bullet^{-} \subseteq \JRad(A_\bullet)$, which will show (i). Since every nilpotent element is left quasi-regular, with $x$ and $y$ as before it suffices to show that $x y$ is nilpotent. Certainly $x$ is nilpotent since it has degree $n < 0$, and since $x\in Z(A_\bullet)$ we have $(xy)^N = x^N y^N = 0$ for $N \gg 0$. This shows that $x \in \JRad(A_\bullet)$.

For any $n \in \bZ$, we define
\[ A_{> n} := \bigoplus_{m > n} A_m, \quad A_{< n} := \bigoplus_{m < n} A_m \]
Let $x \in A_n$ with $n > 0$ and $y = y_{< -n} + y_{-n} + y_{> -n} \in A_\bullet$, where $y_{< -n} \in A_{< -n}$, $y_{> -n} \in A_{> -n}$, and $y_{-n} \in A_{-n}$. Then $x y_{< -n} \in A_\bullet^{-}$, $x y_{> -n} \in A_\bullet^{+}$, and $x y_{-n} \in A_0$. We have shown that $x y_{< -n} \in Z(A_\bullet) \cap \JRad(A_\bullet)$. We claim it suffices to show that $x y_{> -n} + x y_{-n}$ is nilpotent. For if this is the case, then if we set $z_{+} = x y_{> -n} + x y_{-n}$ and $z_{-} = x y_{< -n}$ we have
\[ (z_{+} + z_{-})^N = \sum_{i=0}^N c_{N, i} z_{+}^i z_{-}^{N-i} \]
where $c_{N, i} \in \basefield$ and we have used the fact that $z_{-}$ is central. Since both $z_{+}$ and $z_{-}$ are nilpotent, for $N \gg 0$ each summand will be zero. This shows that $x y$ is nilpotent, hence left quasi-regular and thus in the Jacobson radical. Now we show the claim. Since $y_{-n} \in Z(A_\bullet)$, for any $N > 0$ we have
\[ (x y_{> -n} + x y_{-n})^N = \sum_{i=0}^N M_{N,i}(x y_{> -n}, x) y_{-n}^i \]
where each $M_{N,i}(x y_{> -n}, x)$ is a possibly non-commutative monomial in the elements $x y_{> -n}$ and $x$. Since $x y_{> -n}$ and $x$ are contained in $A_\bullet^{+}$, it follows that $M_{N,i}(x y_{> -n}, x) \in A_{\geq N}$, and hence will be $0$ for $N \gg 0$. This shows that $x y_{> -n} + x y_{-n}$ is nilpotent, proving the claim.
\end{proof}

Note that in \cite[Theorem 4.7]{lauret2003moment}, the structure of graded-semistable Lie algebras is also described. It is more challenging to prove this result via our algebraic approach. As we remarked above, we do not take into account any compatibility with an inner product, which appears to be crucial in the proof of this theorem.

\bibliography{canonical_filtrations_algebras.bib}
\bibliographystyle{alpha}

\end{document}